\documentclass[11pt]{article}
\usepackage{amsfonts}
\usepackage{latexsym,amssymb,amsfonts,graphicx}
\usepackage{amsmath,dsfont}
\usepackage[linewidth=1pt]{mdframed}
\usepackage{mathrsfs}
\usepackage[normalem]{ulem}
\usepackage{epsfig}
\usepackage{tabularx}
\usepackage{geometry}
\usepackage{verbatim}
\usepackage{enumitem}

\usepackage{graphicx}
\usepackage{float}
\usepackage{url}
\usepackage{bm}
\usepackage{color}
\usepackage{natbib}

\usepackage[unicode=true,colorlinks,citecolor=linkcolor,linkcolor=blue,urlcolor=blue]{hyperref}

\providecommand{\U}[1]{\protect\rule{.1in}{.1in}}

\newtheorem{thm}{Theorem}[section]

\newtheorem{lem}{Lemma}[section]

\newtheorem{definition}{Definition}[section]
\usepackage[bf,SL,BF]{subfigure}

\normalsize
\numberwithin{equation}{section}
\allowdisplaybreaks[4]

\definecolor{linkcolor}{rgb}{0,0,0.7}

\begin{document}

\title{Value Maximization under Stochastic Quasi-Hyperbolic Discounting 
}
\author{
Kaixin Yan\thanks{School of Mathematical Sciences, Xiamen University, Xiamen, 361005, China. Email: kaixinyan@stu.xmu.edu.cn}
\and
Wenyuan Wang\thanks{Corresponding author. School of Mathematics and Statistics, Fujian Normal University, Fuzhou, 350007, China; and School of Mathematical Sciences, Xiamen University, Xiamen, 361005, China. Email: wwywang@xmu.edu.cn}
\and
Jinxia Zhu\thanks{Corresponding author. School of Risk and Actuarial Studies, University of New South Wales Sydney, NSW 2052, Australia. Emaill: jinxia.zhu@unsw.edu.au}
}

\date{}

\maketitle
\begin{abstract}
We investigate a value-maximizing   problem incorporating a human behavior pattern: present-biased-ness, for a firm  which navigates strategic decisions encompassing earning retention/payout and  capital injection policies, within the framework of L\'evy processes. We employ the concept of stochastic quasi-hyperbolic   discounting
to capture the present-biased inclinations and model decision making as an intra-personal game with sophisticated decision-makers.  Our analysis yields  closed-form solutions, revealing that double-barrier strategies constitute Markov equilibrium strategies. Our findings reveal that firms, influenced by present-biased-ness,  initiate dividend payments sooner, diminishing overall value compared to scenarios without present-biased-ness (under exponential discounting). We also discuss bailout optimality, providing necessary and sufficient conditions. The impact of behavioral issues is examined in the Brownian motion and jump diffusion cases.

\end{abstract}


\section{Introduction}
In the classical economic theory, when making inter-temporal choices, decision makers  quantify the utility/reward the agent receives in total and use discounting to represent agent's time preference. Generally, a discount function of exponential structure (exponential discounting) is used. However, a phenomenon that has been observed is that  the rate at which
people discount future rewards declines as the length of the delay (waiting time) increases, and equivalently, people's behavior exhibits impatience in the short run and patience in the long run (\cite{Redden(2007)}). Such tendency is called present-biased preferences (\cite{ODonoghueRabin(1999)}). Specifically, 
``individuals are highly impatient about consuming
 between today and tomorrow but are much more patient about
 choices advanced further in the future, for example, between $365$
 and $366$ days from now" (\cite{Barro1999}). 
 Hence, rates of time preference would be
 very high in the short run but much lower in the long run.
In such case, it is more appropriate to use discount functions of hyperbola structure. ``This phenomenon
has been termed hyperbolic discounting by the psychologist Richard Herrnstein" as stated in  \cite{Redden(2007)}.  Most of the existing work on inter-temporal decision making problem uses exponential discounting (generally constant discount rates or constant discount rate in each regime when there are multiple regimes).
Under the classical time preferences model (the exponential discounting model),  choices are time-consistent, that is,
decision makers will make the same utility tradeoff between two periods ($t$ vs. $t+s$) regardless of
when (on or before date $t$) they make the allocation (\cite{Strotz1955}).
Although time consistency has many desirable properties,
it has been criticized for being empirically unrealistic (\cite{CohenEricsonLaibsonWhite2020}). \cite{AinslieHaslam1992} mentioned a directly observed phenomenon: animals and men prefer for ``poorer, earlier alternatives" when they are imminently available. People devalue
the future ``in a curve that is more deeply bowed than  economists' family exponential curve"  (\cite{AinslieHaslam1992}) and research on animal and human behavior concludes that ``discount functions are approximate hyperbolic" (\cite{Laibson1994}).

A discrete-time quasi-hyperbolic discount function has been used to capture the present biasedness and has been used to study a range of behaviors, including consumption, procrastination, addiction,
and job search (\cite{Laibson1997}). 
Under the 
discrete-time quasi-hyperbolic discount model, time is divided into two intervals: the present (consisting of the current period) and the future. Cash-flows in the present period are discounted exponentially with a constant discount factor, $v$. Cash-flows in the future periods are also discounted exponentially with the same discount factor and are further discounted by a constant weight, $\beta$ ($\beta\in(0,1]$). That is, at the reference time point $t=0$, cash-flows in the present period ($t=1$) is discounted by $v$ and those in time $t$ for $t>1$ are discounted by $\beta v^t$.
\cite{HarrisLaibson2013} extended the discrete time discount function to continuous time by introducing a stochastic version of 
quasi-hyperbolic discounting. Under the  stochastic quasi-hyperbolic discounting model, the transition from the present to the future happens with a constant intensity/hazard rate. Specifically, the duration of the present period follows an exponential distribution. Similar to  the discrete-time quasi-hyperbolic discount model, the cash-flows in the present period are discounted exponentially by a constant discount rate, while those in the future period are discounted exponentially by the same discount rate and further multiplied by a constant factor. As the reference time point moves on, the decision maker is modelled as a sequence of selves with different reference points. Each self can only control the decisions/actions for her own present period (starting from her reference point) and do not have control over the decision in her future.

As the first one to  consider alternatives to exponential discounting,
 R. H. Strotz ``recognized that for any discount
function other than exponential, a
person would have time-inconsistent
preferences"
(\cite{FrederickLoewensteinODonoghue2002}), i.e., preferences which imply a conflict between the optimal decisions from the perspectives of the decision maker at different reference time points. Hence, when there is time-inconsistency, it is  not meaningful to look for an optimal policy.  \cite{BjorkMurgoci2010} summarised that there are  different approaches for   handling a family of time inconsistent
problems. One approach is the pre-commitment approach, which is the approach that fix an initial time point and then try to find the optimal control policy that maximize the objective function at the initial point, disregarding the possibilities that at  later points in time, the control law will no longer be optimal. This is like assuming that there is a pre-commitment mechanism to make sure that the decision makers won't change
their mind. 
Another approach is to formulate 
the inter-temporal decision problem as an intra-personal game among selves at different time reference points. The basic idea is that the 
problem is modelled as a non-cooperative game where the decision maker at each time point is viewed as a player and the goal is to find the 
Nash equilibrium strategy for the game.

In this paper, we explore a corporate decision-making problem and investigate the influence of present bias on optimal strategies. Specifically, we consider a  company  that generates  uncertain net earnings at a continuous cost, which is modelled by a spectrally positive L\'evy process.
We investigate the value maximization control problem where the management aims to choose  optimal earning retention/payout and capital injection strategies to maximize the shareholders' value.
We use the discounted cashflow approach and use the expectation of the discounted future dividends as the  profit/shareholders' value. Although there are an extensive amount of works on similar control problems, most of the works in the literature uses exponential discounting 
(\cite{Avanzi11}, \cite{Ya11}, \cite{Bayraktar13}, 
\cite{ZC17}, 
\cite{WW22.b}, \cite{WYZ24}). 
Problems with  present-biased preferences have not received much attention yet in the literature except for a few works on applying non-exponential discounting to pure diffusion models or the compound Poisson model. 
All the works taking into account of the  present biased preferences of the decision makers are based on either compound Poisson (\cite{Li16}, \cite{CW16}), 
or in one reference, Brownian motion (\cite{ZW14}), 
or diffusion processes (\cite{CL18}, \cite{ZSY2020}).
We will investigate the stochastic control problem for a general L\'evy model which includes the compound Poisson and Brownian motion models as special cases. L\'evy processes are stochastic processes characterized by their remarkable mathematical properties, including independent and stationary increments, and have been widely applied in various disciplines due to their ability to capture complex and non-Gaussian behaviors. In physics, L\'evy processes are crucial to the study of turbulence, laser cooling and in quantum field theory. In engineering, L\'evy processes play an important role in the investigation  of networks, queues and dams. 
In the field of mathematical finance, L\'evy processes, including their special cases, are extensively employed to provide a more accurate description of financial markets compared to Brownian motion. Therefore, it is imperative to conduct a study on the impact of present bias on dividend payment decisions based on the Lévy process. 

In this work, we base our study on  spectrally positive L\'evy  processes with negative drifts, which are L\'evy processes with positive jumps and negative drifts. These are appropriate models for  situations,  when  a company is driven by inventions or discoveries (\cite{Bayraktar13}, 
\cite{ZC17}). For these companies, continuous expenditures are made to sustain operations in exchange for random earnings.   We consider a decision maker who aims to maximize the company's value, which is measured by the expected present value of all the future cashflows. We follow \cite{HarrisLaibson2013} by using a particular form of non-exponential discounting, the quasi-hyperbolic stochastic discounting, to formulate the decision makers' present biased-ness. We assume  that  the decision makers are sophisticated and that there is no commitment ability.
Assuming sophisticated decision-makers, we formulate the optimization problem as an intra-personal game among the decision markers themselves at different reference points and the goal is to seek stationary Markov equilibrium strategies for  this intra-personal game. We solve the intra-personal game problem, and find closed-form solutions for Markov equilibrium strategies within the  L\'evy model. We reveal that that a distinct set of strategies, known as double-barrier strategies, constitutes a Markov equilibrium strategy. We   examine the impact of present-biasedness on decisions, focusing on understanding  the deviations of these equilibrium strategies from those made by rational agents and the potential impacts on the surplus dynamics and the firm values. 
 Our findings reveal that decisions under present bias exhibit impatience in paying dividends, prompting firms to initiate dividend payments earlier (reflected in the lower dividend barrier) compared to their exponential discounting counterparts, consequently diminishing the firm's overall value. Additionally, we delve into the optimality of bailouts, providing a necessary and sufficient condition for optimality. 
Further discussion extends to the detailed examination of the impact of behavioral issues in two specific examples: the Brownian motion and jump diffusion cases.  This exploration offers valuable insights into the impact of behavioral issues on decisions and firm values.

In this paper, we systematically examine the dynamics of equilibrium cashflow management strategies within the framework of stochastic quasi-hyperbolic discounting. A key insight emerges from the examination of the equilibrium dividend threshold concerning varying degrees of present bias: decision-makers exhibit increased impatience, leading to a desire for earlier dividend payments as the discounting for future cashflows deviates from exponential discounting or when future arrivals become more likely. Additionally, our investigation explores the impact of behavioral issues on incurred losses, revealing a consistent trend of losses attributed to behavior issues, with the magnitude increasing as present-bias intensifies. To provide a comprehensive understanding, we simulate paths for the uncontrolled L\'evy process and the corresponding optimally controlled processes under both stochastic quasi-hyperbolic and exponential discounting. Our exploration of behavioral issues sheds light on the consistent presence of losses, emphasizing the significant role of present-bias in influencing firm values. Overall, this study contributes to the understanding of optimal financial decision-making under stochastic quasi-hyperbolic discounting, providing insights for strategic financial management in real-world scenarios.

The layout of the rest of the paper is as follows. In Section \ref{sec:model1}, we present the definition of L\'evy processes and some of its characteristic properties. We introduce the dynamics of the stochastic process under consideration and provide the definition of stochastic quasi-hyperbolic discounting. Furthermore, we formulate the stochastic control problem as an intra-personal non-cooperative game, with the objective of seeking Markov equilibrium strategies. Moving forward, in Subsection \ref{auxilary}, we examine a special class of candidate strategies - the double-barrier strategies,  and solve the  stochastic control problem in Subsection \ref{Markov}.  Our findings reveal that a double-barrier strategy, with strategically placed barriers, represents a desired Markov equilibrium strategy within the general L\'evy model.  We discuss the impact of present-biased-ness on the firm value in Subsection \ref{subsect3.3}. We explore the optimality of bailouts and provide necessary and sufficient conditions for such optimality in Section \ref{contionforcapitalinjections}.
In Sections \ref{examples} and \ref{jumpdiffusion}, we consider two classical models, which are special cases of the general L\'evy model.
We utilize the analytical results obtained for the general setting to analyze these specific cases. Additionally, we offer numerical illustrations to thoroughly explore the implications of present-bias and derive intriguing insights. We conclude the paper with some comments and remarks in Section \ref{conclusion}.

\section{Problem Formulation }\label{sec:model1}
Consider a  company  that makes continuous expenditures/investment, generating   uncertain (random)  earnings. Let $X_t$ represent the cumulative net earnings by time $t$.
The company chooses policies on retaining earnings and distributing dividends, and the retained earnings stay in  the cash reserve/ surplus. We use $L{(t)}$ to represent the cumulative amount of dividends paid out up to time $t$. Here, if there is a lump sum payment at time $t$, we assume that the amount is deducted immediately before time $t$, that is, any jump of $L$ at time $t$ is assumed to happen from $t-$ to $t$ and the process $L:=\{L{(t)}; t\ge 0\}$ is right continuous. Furthermore, it is natural to assume that all dividend payments need to be non-negative and hence, $L$ is non-decreasing.  The company also controls capital injections (for example, via equity issuance). Capital injections are costly and for every unit of capital raised, the total cost is $\phi$ ($\phi>1$).   We use $R{(t)}$ to represent the cumulative amount of capitals injected by time $t$. 
The amount of each capital injection must be positive and we assume that $R:=\{R{(t)};t\ge 0\}$ is non-decreasing and right continuous with $R({0-})=0$. Furthermore, we consider a bailout problem where the beneficiaries of dividends are supposed to inject capitals into the surplus process so the resulting surplus processes is always non-negative. We use $U^{R,L}$ to represent the cash reserve (the controlled reserve/surplus process) under the strategy $(R,L)$:
\begin{align}
    U^{R,L}{(t)}:=X{(t)}-L{(t)}+R{(t)}, \quad t\geq 0.
    \end{align}

The company wants to choose policies so that it's profit/shareholders' value is maximized. Specifically, the value of a policy is quantified by the expected total discounted dividends net of the total cost of capital injections. For $x\in\mathbb{R}$, we denote by $\mathrm{P}_x$ the law of $X$ starting from $X(0)=x$ and write $\mathrm{E}_x$ the associated expectation. 
In order to value the cashflows, we define $P(t,s)$ to be the present value at time $t$ of 1 unit of cash at time $s$.
So the expected discounted value at time $t$ of all the current and future dividends  net of the discounted value of the capital injections plus the associated cost   is
\begin{align*}
\mathrm{E}_{x}\left[\int_{t}^{\infty}P(t,s)\left({d}L{(s)}-\phi{d}R{(s)}\right)\right].
\end{align*}

We consider the situation where decision makers are present-biased and follow the quasi-hyperbolic stochastic discounting framework in \cite{HarrisLaibson2013}. For any time point $t$, the present period of the decision maker lasts for a duration of $\tau$, which follows an exponential distribution with mean $1/\lambda$. During the present period, the discount factor to the decision maker at time $t$ for any cashflow at time $s$ ($s$ is in the present period) is $e^{- \delta(s-t)}$, while outside the present period, the cashflow at time $s$ (for $s>t+\tau$) will be discounted by the time-$t$ decision maker by a further discount factor $\beta$ multiplied by the usual discount factor $e^{- \delta(s-t)}$. That is,
\begin{align}
P(t,s)=
\begin{cases}e^{- \delta(s-t)},& t<s<t+\tau,\\
\beta e^{- \delta(s-t)},& s>t+\tau.\end{cases}
\end{align}
The above $(\beta,\delta)$ formulation, with $\beta\in(0,1]$ and $\delta\geq0$, captures some essence and many qualitative implications. It nests exponential discounting: if $\lambda=0$, the future never arrives; if $\beta=1$, there is no difference in discounting between present and future (\cite{HarrisLaibson2013}). 

Under the above stochastic formulation for discounting, the control problem can be considered as  a problem with a sequence of decision makers (autonomous selfs), each of which makes decisions during her own present period, and cares about the total profit in the present and future but does not control the policies in the future. Such a problem  leads to an intra-personal game and we follow the literature to employ the Markov-perfect equilibrium concept (\cite{HarrisLaibson2013}). 
 Markov-perfect equilibrium (MPE) strategies are feedback
control laws that maintain Markov structure.  A Markov strategy is a strategy that  is measurable with respect to the state space.
That is, the control at time $t$ depends on time $t$ and the state at time $t$ only: {$dL{(t)}$ and $dR{(t)}$ can be written as (deterministic and measurable) functions of $(t,U^{R,L}{(t-)})$}. 
Intuitively, Makov strategies only ``depend on information that is directly payoff relevant", which, in this case, is the current level of cash reserve. MPE is a refinement of subgame-perfect equilibrium   in which all players use Markov strategies.
MPE is a refined version of Nash equilibrium in game theory.  As commented in \cite{MaskinTirole2001}, MPE embodies both practical and theoretical  value.   It simplifies dynamic game models and aligns with bounded rationality through its use of simple Markov strategies. MPE uniquely captures the diminishing importance of past actions, and adheres to the principle that only significant past elements should strongly impact behavior.
We further restricting the strategies to be  stationary
MPE (SMPE). That is, we consider strategies all
selves use the same strategy.

We use $((R,L),(\tilde{R},\tilde{L}))$ to represent the strategy when the current self  adopts $(R,L)$ and all the future selves employ $(\tilde{R},\tilde{L})$.  Let $D^{((R,L),(\tilde{R},\tilde{L}))}{(t)}$ represent the cumulative amount of withdrawals from up to $t$ under the strategy $((R,L),(\tilde{R},\tilde{L}))$. Then,   $d D^{((R,L),(\tilde{R},\tilde{L}))}{(t)}=d L{(t)}$ for $t\in(0,\tau)$ and  $d D^{((R,L),(\tilde{R},\tilde{L}))}{(t)}=d \tilde{L}{(t)}$ for $t\in[\tau,\infty)$.
The reward to current self  is  the  expected present value at time $0$ of all the future net consumption/dividends received. 
Given $U({0-})=x$, for any $x\ge 0$, the current self's expected reward function of the strategy $((R,L),(\tilde{R},\tilde{L}))$ is
\begin{align}
&\mathcal{P}(x;(R,L),(\tilde{R},\tilde{L}))\nonumber\\
=&\mathrm{E}_{x}\bigg[\int^{\tau}_{0} e^{-\delta t}(d L{(t)}-\phi dR{(t)})
+
\beta \int_{\tau}^{\infty} e^{-\delta  t}(d\tilde{L}{(t)}-\phi d\tilde{R}{(t)})\bigg].\label{PLb}
\end{align}
Note that $\tau$ in the above equation, \eqref{PLb}, is the time at which control passes from the current self at time $0$ to the next self. The first integral inside the expectation notation sums all the discounted dividends net of the capital injections received during the present period under the consistent discounting at the rate of $\delta$, and the second integral collects all the discounted value at time $0$ of dividends minus the capital injections from time $\tau$. If we use the  double expectation formulae and the strong Markov property, the expected value of the second term inside the expectation notation in \eqref{PLb} can be rewritten as
\begin{align}
&\mathrm{E}_{x}\left[
\int_{\tau}^{\infty} e^{-\delta  t}(d\tilde{L}{(t)}-\phi d\tilde{R}{(t)})\right]\nonumber\\
=&\mathrm{E}_{x}\left[e^{-\delta \tau}\left.
\mathrm{E}_x\left[\int_{\tau}^{\infty} e^{-\delta  (t-\tau)}(d\tilde{L}{(t)}-\phi d\tilde{R}{(t)})\right|\mathcal{F}_{\tau}\right]\right]\nonumber\\
=& \mathrm{E}_{x}\bigg[e^{-\delta \tau} \mathrm{E}_{U^{R,L}{(\tau)}}
\bigg[\int_{0}^{\infty} e^{-\delta  t}(d\tilde{L}{(t)}-\phi d\tilde{R}{(t)})\bigg]\bigg].\label{PLb-2}
\end{align}
Define   $\mathcal{P}^E(\cdot;(\tilde{R},\tilde{L}))$ as the expected payoff function under exponential discounting with discount rate $\delta$:  \begin{align}
&\mathcal{P}^E(x;(\tilde{R},\tilde{L}))=\mathrm{E}_{x}\left[\int_{0}^{\infty} e^{-\delta  t}(d\tilde{L}{(t)}-\phi d\tilde{R}{(t)})\right]. \label{14116-3}\end{align}
Note that $\tau$ is an independent exponential random variable. By applying the strong Markov property of $(\tilde{R},\tilde{L})$ and  combining \eqref{PLb}, \eqref{PLb-2}  and \eqref{14116-3} we can obtain 
\begin{align}
&\mathcal{P}(x;(R,L),(\tilde{R},\tilde{L}))
=\mathrm{E}_{x}\bigg[\int^{\tau}_{0} e^{-\delta t}(d L{(t)}-\phi dR{(t)})
+
\beta e^{-\delta \tau} \mathcal{P}^E(U^{R,L}{(\tau)};(\tilde{L},\tilde{R}))
\bigg].\label{PLb-3}
\end{align}

As noted earlier  that under every non-exponential discounting, the control problem is  time inconsistent
problem and  the Bellman optimality principle does not hold. Under the quasi-hyperbolic stochastic discounting,  the decision maker  is modeled as a sequence of
autonomous selves. Each self controls capital injection and dividend distribution actions  during her own present period only and does not control the actions in her future periods, although she cares about the future decisions. This formulation is an intrapersonal game and we will seek stationary Markov-perfect equilibrium by following  the literature in
this area, as it is not meaningful to seek an optimal solution, which is hard to define in such situation. 
Under such formulation, the  players (the selves) make decisions sequentially, and each player's (self's) actions can influence the evolution of the game in future periods. Each player's strategy specifies their actions at each possible state of the game during their present period taking into account the player's beliefs about how the future selves will behave in the future states. The MPE strategies chosen by the players are optimal (in term of maximizing the total value, including the value generated in the current and future periods) given their beliefs and the current state of the game, accounting for the potential consequences of their actions on the evolution of the game.

A strategy $\pi=(R,L)$ is called \textit{admissible} if both $L=\{L(t);t\ge 0\}$ and $R=\{R(t);t\ge 0\}$
are positive increasing cadl\'ag processes, adapted to  $(\mathcal{F}_t)_{t\ge 0}$,  the filtration generated by $\{X(t);t\ge 0\}$, and $X^\pi(t)\ge 0$ for all $t\ge 0$. 
Let $\Pi$ represent the set of all the admissible strategies.

For any decision maker in the present period, the objective is to find an admissible stationary Markov-perfect equilibrium (MPE) strategy, $\pi^{\ast}=({R}{^\ast},L{^\ast})$, such that $$\mathcal{P}(x; (R^{\ast},L^{\ast}),(R^{\ast},L^{\ast}))=\sup_{(R,L)\in\Pi}\mathcal{P}(x;(R,L),({R}{^\ast},L{^\ast})).$$ 
The strategy $\pi^\ast$ is optimal in the sense that, for all $x$, $\pi^*$ maximizes $ \mathcal{P}(x; \pi,(R^{\ast},L^{\ast}))$ with respect to $\pi$ in the set of admissible strategies.
This strategy is the one under which  the current self achieve its objective to maximizes her value  anticipating that the future selves intend to maximize their future values from their perspective.  
If an admissible stationary MPE (denoted by $(R^{\ast}, L^{\ast})$) exists, we define the value function $V(x)=\mathcal{P}(x; (R^{\ast},L^{\ast}),(R^{\ast},L^{\ast}))$.

We assume that stochastic process $X=\{X{(t)}; t\geq 0\}$ is a L\'{e}vy process defined on the probability space $(\Omega, \mathcal{F},\mathrm{P})$. We assume that $X$ is a $(-\gamma,\sigma,v)$-L\'evy process, where $\gamma$ is a real number, $\sigma$ is non-negative, and $v$  is a measure concentrated on $\mathbb{R}^+ \symbol{92} \{0\}$ such that 
$$\mathrm{E}[e^{iqX{(t)}}]=\exp\left({-t \left(-i\gamma q +\frac12\sigma^2 q^2 +\int_{\mathbb{R}^{+}\setminus \{0\}}(1-e^{iqy}+iqy  \mathbf{1}_{\{y<1\}})v({d}y)\right)}\right).$$ The process $X$ is a \textit{spectrally positive}  L\'{e}vy process as it has no negative jumps. 
As in the standard literature, we assume \begin{align*}
\int_0^{\infty}(1\wedge y^2)\upsilon({d}y)<\infty.
\end{align*}

As pointed out by Bernt Øksendal in the  lecture notes for An Introduction to Stochastic Control, with
Applications to Mathematical Finance
(\cite{Oksendal2015}) 
that the above condition still ``allows for many interesting kinds of L\'evy process" including the one that has infinitely many small jumps, which is the L\'evy process satisfying
$\int_{\mathbb{R}^{+}\setminus \{0\}} (1\vee y)v(dy)=\infty.$  L\'evy processes are continuous time processes that are right continuous and have left limits. Furthermore, L\'evy processes have independent and stationary increments (\cite{Bertoin96}). 
 These processes can capture the characteristics of asset price processes effectively. One well-known example is the Black-Scholes model, introduced by \cite{BS73}, which is based on a L\'evy process. The model provided a solid mathematical foundation for options trading, revolutionizing finance theory and practice. \cite{Me76} proposed a more general version of L\'evy processes known as the Merton model. The discontinuous price process became essential for understanding how firms determine their capital structure. Another notable L\'evy process model is the jump-diffusion model, introduced by \cite{Ko02}. The model incorporates jumps, similar to Merton's, with the jump size being double-exponentially distributed. It has been widely used in various option-pricing problems. 
In addition to these models, there are several others popular L\'evy processes-based models in the mathematical finance. These include the generalized hyperbolic model and its special case, the normal inverse Gaussian distribution, the CGMY L\'evy process, as well as the Meixner process, etc. (\cite{Pa08}). 

According to Chapter I in \cite{Bertoin96}, the L\'evy process $X$ can be represented in the following form:
\begin{align}
    X(t)=\gamma t+\sigma B(t) + X_1(t)+M_1(t),\label{decompostion}
 \end{align}
where  
$B(t)$ represents a standard Brownian motion, $X_1(t)$ is a compound Poisson process of jumps of magnitude greater than unity, and $M_1(t)$ is a
square integrable martingale with an almost surely countable number of jumps
on each finite time interval which are of magnitude less than unity. The three stochastic processes, $B$, $X_1$ and $M_1$, are mutually independent. 
Note that the generator  of the L\'evy process $X$, denoted by $\mathcal{A}$, is:
\begin{align}
    (\mathcal{A}f)(x)
    &=\frac{\sigma^2}2 f^{\prime\prime}(x)+\gamma f^\prime(x)+\int_0^\infty \left(f(x+y)-f(x)-f^\prime(x)y\mathbf{1}_{\{y<1\}}\right)v(dy).\label{generator}
\end{align}


It is well-known that  $X$ has paths of bounded variation if and only if $\sigma=0$ and  $\int_{0}^{1}z\upsilon({d}z)<\infty$.
As convention, we rule out the case that $X$ has monotone paths (i.e., $X$ is a subordinator), and so we assume that when $X$ is of bounded variation: \begin{align}
-\gamma+\int_{0}^{1}z\upsilon({d}z)>0.
\end{align}
To exclude the trivial case, we further assume throughout the paper that
\begin{align*}
\mathrm{E}[X(1)]
<
\infty.
\end{align*}

\section{Equilibrium results} \label{section3}
Note that  capital injections need to be made to ensure that the surplus does not become negative and  capital injections are costly. 
It is natural to  conjecture that it might be optimal not to inject capitals unless it is  absolutely necessary and the amount of each capital injection should be kept at  the minimal amount. That is, the company injects capital only when the surplus process reaches $0$ from above and the amount of injection is minimal so that it keeps pushing the L\'evy process to stay at or above $0$.  For the dividend payout, we will first consider a smaller subset of  strategies that has been shown to contain the optimal strategies in many settings and explore whether it contains MPE.

\subsection{A class of double-barrier strategies}
\label{auxilary}
We call the strategies defined below double-barrier strategies. 

\begin{definition}(Double-barrier strategies) Let $\pi^{0,b}=(R^0, L^b)$ represent the dividend and capital injection strategy such that when the controlled surplus process
 attempts to exceed the barrier $b$, the excess part beyond the barrier  is paid out as dividends and the process then remains ``stuck" at the barrier until it changes direction and moves downward away from the barrier $b$, and that when the controlled process
 attempts to down-cross  the barrier $0$, capital injections equal to the deficit part below $0$ are provided and the process then remains ``stuck" at the barrier $0$ until it changes direction and moves upward.
 \end{definition}
 


For convenience, we use $V_{0,b}$ to represent the expected payoff under the situation that the current self chooses the double barrier strategy $\pi^{0,b}$ and the future selves use the same strategy $\pi^{0.b}$: 
\begin{align}\label{expdouble}
V_{0,b}(x):=&\mathcal{P}(x;\pi^{0,b},\pi^{0,b})\\
  =&\mathrm{E}_{x}\bigg[\int^{\tau}_{0} e^{-\delta t}(d L^b{(t)}-\phi dR^0{(t)})
+\beta e^{-\delta \tau} \mathcal{P}^E(U^{\pi^{0,b}}{(\tau)};\pi^{0,b})
\bigg]\nonumber\\
=&\mathrm{E}_{x}\bigg[\int^{\infty}_{0} e^{-(\delta+\lambda) t}(d L^b{(t)}-\phi dR^0{(t)})
+
\lambda \beta \int^{\infty}_{0} e^{-(\delta +\lambda) t} \mathcal{P}^E(U^{\pi^{0,b}}{(t)};\pi^{0,b})dt
\bigg].\label{def2-v0b}
\end{align}
We further define
\begin{align}
V_{0,b}^E(x):=\mathcal{P}^E(x;\pi^{0,b}).\label{def-vobE}
\end{align}
The function, $V_{0,b}^E(x)$, measures the performance of the double barrier strategy $\pi^{0,b}$ under exponential discounting with the discount rate $\delta$. This function is same as the expected payoff under the same strategy in \cite{Bayraktar13}. 
We will present some properties satisfied by $V_{0,b}^E(x)$, which will be used in deriving the main results.  We start with defining some standard quantities of L\'evy processes and in particularly, the scale functions and recalling some of their associated properties. Let  $\psi: [0,\infty)\rightarrow \mathbb{R}$ denote the  Laplace exponent of the  L\'evy process $X$, then $\psi$ satisfies
\begin{align*}
\mathrm{E}[e^{-\theta X{(t)}}]=:e^{\psi(\theta)t}, \quad t,\theta\geq 0,
\end{align*}
and according to the L\'{e}vy-Khintchine formula
\begin{align*}
\psi(\theta)=-\gamma \theta+\frac{\sigma^2}{2}\theta^2+\int_{(0,\infty)}(e^{-\theta z}-1+\theta z\mathbf{1}_{\{|z|<1\}})\upsilon(dz), \quad \theta\geq 0.
\end{align*}
Following are some standard results  in the literature (see, for example, \cite{Bayraktar13}).
Let us also recall the $q$-scale function for the spectrally positive L\'{e}vy process $X$. For any $q>0$, there exists a continuous and increasing function $W_{q}:\mathbb{R}\rightarrow [0,\infty)$ such that \begin{align}
    W_{q}(x)=0,\quad x<0,\label{14823-7}  \end{align} and its Laplace transform on $[0,\infty)$ is given by
\begin{align*}
\int_0^{\infty} e^{-sx}W_{q}(x)dx=\frac{1}{\psi(s)-q},\quad s>\Phi(q),
\end{align*}
where $\Phi(q):=\sup\{s\geq 0: \psi(s)=q\}$. The function $W_{q}$ is called the $q$-scale function. Furthermore, we can construct two families of functions, $Z_{q}(x)$ and $\overline{Z}_{q}(x)$,  by
\begin{align}
&Z_{q}(x):=1+q\int_0^x W_{q}(y)dy,\quad x\in\mathbb{R},\label{8823-3}\\
&\overline{Z}_{q}(x):=\int_0^x Z_{q}(y)dy,\quad x\in\mathbb{R}.\nonumber
\end{align}
Notice that because $W_q$ is $0$ on the negative half line and so 
\begin{align}
    Z_q(x)=1\quad \mbox{and}\quad \overline{Z}_q(x)=x,\quad x\le 0.\label{9823-9}
\end{align}

If $X$ is of paths of bounded variation, $W_{q}(x)\in C^1(0,\infty)$ if and only if the L\'{e}vy measure $\upsilon$ has no atoms.
Moreover, if $\sigma>0$, we have $W_{q}(x)\in C^2(0,\infty)$.
We also know that
\begin{align}
\begin{split}
W_{q} (0) &= \left\{ \begin{array}{ll} 0 & \textrm{if $X$ is of unbounded
variation,} \label{9823-4}\\ \frac1{-\gamma+\int_{(0,1)}z\upsilon(dz)}
& \textrm{if $X$ is of bounded variation,}
\end{array} \right.
\end{split}
\end{align}
and 
\begin{align}
W_{q}^\prime (0+) &:=\lim_{x\downarrow 0}W_{q}^\prime (x)=
\begin{cases}
\frac{2}{\sigma^2} & \mbox{if } \sigma>0,\\
\infty&\mbox{if } \sigma=0 \mbox{ and } v(0,\infty)=\infty,\\
\frac{q+v(0,\infty)}{\left(-\gamma+\int_{(0,1)}z\upsilon(dz)\right)^2}
&\mbox{if $X$ is compound Poisson}.    \end{cases}
\end{align}

Define
\begin{eqnarray}\label{bE.def}
b^E:=\inf\{b>0:Z_{\delta}(b)-\phi\geq0\}.
\end{eqnarray}
Notice that the function $Z_{\delta}(x)$ is strictly increasing in $x$ on $(0,\infty)$, that $Z_{\delta}(0)=1<\phi$, and that $\lim\limits_{x\rightarrow\infty}Z_{\delta}(x)=\infty$. So $b^E\in(0,\infty)$, and thus,  $b^E$ is well defined and satisfies
\begin{align}
Z_{\delta}(b^E)-\phi=0.\label{9823-5}
    \end{align}

As noted earlier that $V_{0,b}^E(x)$ is the same as the function $\bar{v}_b$ in \cite{Bayraktar13}. Furthermore, we can see that the barrier $b^E$ defined above is same as the $b{^\ast}$ defined in  \cite{Bayraktar13}, and  $V_{0,b^E}^E(x)$ is same as $\bar{v}_{b{^\ast}}$ in  \cite{Bayraktar13}. The following lemma collects some results that have been proved in   \cite{Bayraktar13}.
\begin{lem}\label{Vb.E} 
(i) For any $b\ge 0$, the function $V_{0,b}^E(x)$ has the following representation \begin{eqnarray}\label{V.0bE.x}
V^E_{0,b}(x)
\hspace{-0.3cm}&=&\hspace{-0.3cm}
\left\{
\small
\begin{aligned}
&
-
\overline{Z}_{\delta}(b-x)
-\frac{\psi^{\prime}(0+)}{\delta}
+\frac{Z_{\delta}(b-x)}{\delta W_{\delta}(b)}
\Big[Z_{\delta}(b)-\phi
\Big],& x\in[0,b],\\
&x-b+V^E_{0,b}(b),& x\in(b,\infty),\\
&\phi x+V^E_{0,b}(0), &x\in(-\infty,0).
\end{aligned}
\right.
\end{eqnarray}

\noindent (ii)
For any $b\in(0,\infty)$, the function $V_{0,b}^E(x)$ is strictly increasing and continuously differentiable on $(-\infty,\infty)\setminus\{b\}$. And, if $X$ has paths of unbounded variation, $V_{0,b}^E(x)$ is continuously differentiable on  $(-\infty,\infty)$ and twice continuously differentiable over $(0,\infty)\setminus\{b\}$. Furthermore, we have 
\begin{align}\label{VE.hjb}
(\mathcal{A}-\delta)(V^{E}_{0,b})(x)=0,\quad x\in(0,b).
\end{align}
(ii)
The function $V_{0,b^E}^E(x)$ is continuously differentiable on $(-\infty,\infty)$ and is concave over $(0,b^E)$ with
\begin{align}
1<V_{0,b^E}^{E\,\prime}(x)<\phi,\quad x\in(0,b^E),\\
V_{0,b^E}^{E\,\prime}(x)=1, \quad x\in[b^E,\infty),\label{9823-7}\\
V_{0,b^E}^{E\,\prime}(x)=\phi, \quad x\in(-\infty,0).
\end{align}
Furthermore, if $X$ has paths of unbounded variation, $V_{0,b^{E}}^E(x)$ is twice continuously differentiable on $(0,\infty)$.
\end{lem}

\begin{lem}\label{lem2.4}
(i) For any $b\in[0,b^E]$, the function $V_{0,b}^{E}(x)$ is concave on $[0,b)$.\\
\noindent (ii) If $X$ is of unbounded variation, $V_{0,b}^{E\,\prime\prime}(b^E-)=0$. 
\end{lem}

We can now derive expressions for $V_{0,b}$, the value under stochastic quasi-hyperbolic discounting when the current and all the future selves employ the double-barrier strategy, $\pi^{0,b}$, in terms of the scale functions as well.

\begin{lem}\label{Vb}
Fix $b\in(0,\infty)$. The function $V_{0,b}(x)$ has the following representation
\begin{eqnarray}
\label{v0b(x)-1}
\hspace{-0.3cm}
V_{0,b}(x)
\hspace{-0.3cm}&=&\hspace{-0.3cm}
\left\{
\small
\begin{aligned}
&
-
\overline{Z}_{\delta+\lambda}(b-x)
-\frac{\psi^{\prime}(0+)}{\delta+\lambda}
+
\frac{\lambda\beta}{\delta+\lambda}\bigg[V_{0,b}^E(0)+
\int_{0}^{b}V_{0,b}^{E\,\prime}(y)
Z_{\delta+\lambda}(y-x){d}y\bigg]&
\\&
+\frac{Z_{\delta+\lambda}(b-x)}{(\delta+\lambda)W_{\delta+\lambda}(b)}
\Big[Z_{\delta+\lambda}(b)-\phi
-\lambda\beta\int_{0}^{b}V_{0,b}^{E\,\prime}(y)W_{\delta+\lambda}(y)
{d}y\Big],& x\in[0,b],\\
&x-b+V_{0,b}(b),& x\in(b,\infty),\\
&\phi x+V_{0,b}(0), &x\in(-\infty,0).
\end{aligned}
\right.
\end{eqnarray}
The function $V_{0,b}(x)$
is continuously differentiable on $(-\infty,\infty)\setminus\{b\}$.
And, if $X$ has paths of unbounded variation, $V_{0,b}(x)$ is once continuously differentiable on $(-\infty,\infty)$ and twice continuously differentiable on $(0,\infty)\setminus\{b\}$ with
\begin{align}
    V_{0,b}^\prime(b)=1.\label{9823-8} 
\end{align} Furthermore, we have
\begin{eqnarray}\label{V.hjb-1}
\mathcal{A}V_{0,b}(x)-(\delta+\lambda) V_{0,b}(x)+\lambda\beta V_{0,b}^E(x)=0,\quad x\in(0,b).
\end{eqnarray}
\end{lem}

The proof of the above lemma is provided in the Appendix. We establish the explicit expression  
 by leveraging some standard results on L\'evy process from the literature. It involves deconstructing an expectation into separate integrals, incorporating known results from literature, and employing integration by parts to simplify the expression. The derived formula is then compared with existing findings, leading to the establishment of a system of differential equations. 

The formulae, \eqref{V.hjb-1}, in the last Lemma can be rewritten as
\begin{eqnarray}\label{V.hjb}
\mathcal{A}V_{0,b}(x)-\delta V_{0,b}(x)+\lambda(\beta V_{0,b}^E(x)-V_{0,b}(x))=0,\quad x\in(0,b).
\end{eqnarray}
Recall that $\lambda$
represents the instantaneous intensity/hazard rate of making the transition from the ``present" 
to the ``future". As noted by \cite{HarrisLaibson2013}, at the transition point,  the continuation value equivalent to $V_{0,b}^E(x)$  begins and this is discounted by a discount function that is a fraction, $\beta$, of  the discount function in the present. If $\beta=1$ , then there is no difference in the discounting for the cashflows in the present and future periods. If $\lambda=0$, there are no future periods. 

\subsection{Markov equilibrium strategies}\label{Markov}
Now we will define a quantity, $b{^\ast}$, which will be shown to be the optimal upper barrier in the sense that the double-barrier strategy with $b{^\ast}$ as its upper barrier (for both the current and future selves) dominates, in terms of the expected total value, any other double-barrier strategies with $b\neq b{^\ast}$ as its upper barrier (for both the current and future selves). We will then investigate whether the optimal double-barrier strategy is a Markov equilibrium strategy. 
\begin{definition}
    Define
\begin{eqnarray}
\label{smoothcond.1}
b^{\ast}:=\inf\{b>0:Z_{\delta+\lambda}(b)-\phi
-\lambda\beta\int_{0}^{b}V_{0,b}^{E\,\prime}(y)W_{\delta+\lambda}(y)
{d}y\geq0\}.
\end{eqnarray}
By convention, $b^\ast =\infty$, if $Z_{\delta+\lambda}(b)-\phi
-\lambda\beta\int_{0}^{b}V_{0,b}^{E\,\prime}(y)W_{\delta+\lambda}(y)
{d}y<0$ for all $b>0$.
\end{definition}

Through a lengthy  derivation (see Appendix), we can show $Z_{\delta+\lambda}(b^{E})-\phi-\lambda\beta\int_{0}^{b^{E}}V_{0,b^{E}}^{E\,\prime}(y)W_{\delta+\lambda}(y){d}y\geq0$, which implies  that  $b{^\ast}$ bounded by $b^E$, the optimal upper barrier in the exponential discounting counterpart.
\begin{lem}\label{lem.3.5}
 We have $0<b^{\ast}\le b^E<\infty$.
\end{lem}

From above we know that $b^*$ is attainable and hence it is the smallest positive root of $$\ell(x):=Z_{\delta+\lambda}(x)-\phi
-\lambda\beta
\int_{0}^{x}V_{0,x}^{E\,\prime}(y)
W_{\delta+\lambda}(y){d}y,\quad x\in[0,b^E].$$

Through an extensive proof (see Appendix) by integrating theoretical results and mathematical techniques, we can establish the following smoothness of the value function under the double-barrier strategy with dividend barrier $b^\ast$ and prove the concavity of 
the function by showing that its derivative function is non-increasing.

\begin{lem}\label{b.ast}
(a) The function $V_{0,b^{\ast}}(x)$ is continuously differentiable on $(-\infty,\infty)$ and is concave on $(0,\infty)$.
(b) Furthermore, if $X$ has paths of unbounded variation, $V_{0,b^{\ast}}(x)$ is twice continuously differentiable on $(0,\infty)$. (c) If $X$ has paths of bounded variation, $V_{0,b^{\ast}}^\prime(b^\ast)=1$.
\end{lem}

Noting the smoothness of the function $V_{0,b^{\ast}}$, by applying the generalized It\^o formulae (see, Theorem 4.57 in \cite{JaSh2003}) to $e^{-(\delta+\lambda)t}V_{0,b^{\ast}}(U^{R,L}{(t)})$ for any admissible strategy $(R,L)$ and utilizing the properties and results obtained earlier, through a lengthy and complicated process we can show that $V_{0,b^{\ast}}(x) \ge \mathcal{P}(x,\pi,\pi^{0,b^{\ast}})$ for any admissible strategy $(R,L)$. This implies that the function $V_{0,b^{\ast}}(x)$ dominates the supremum of the performance function over all the admissible strategies. Since $V_{0,b^{\ast}}(x)$ is the performance function associated with a particular admissible strategy, $\pi^{0,b^*},$ we can conclude that  $\pi^{0,b^\ast}$ is an MPE strategy.

\begin{thm}\label{16923-1}
The double-barrier dividend and capital injection strategy $\pi^{0,b^{\ast}}$ is an MPE strategy, that is, $V_{0,b^{\ast}}(x)=\mathcal{P}(x;\pi^{0,b^{\ast}},\pi^{0,b^{\ast}})=\sup_{\pi\in {\Pi}}\mathcal{P}(x,\pi,\pi^{0,b^{\ast}})$.
\end{thm}

The above theorem suggests that the double-barrier strategy, $\pi^{0,b^\ast}$, which prescribes to inject capitals at the minimal rate when and only when the current reserves touches $0$ so as to prevent the reserve/surplus falling below $0$ and to pay out dividends when and only when the reserve/surplus reaches or exceeds $b^\ast$. The controlled stochastic process under the strategy $\pi^{0,b^\ast}$ is a double reflected L\'evy process that  is constrained to stay within $0$ and $b^\ast$. 

 As noted earlier,  $b^\ast \le b^E$.  This implies that under the equilibrium strategy  $\pi^{0,b^{\ast}}$, the threshold for paying dividends in the present-biased case is lower than the exponential discounting case. Consequently,  firms will start  dividend payments earlier under the quasi-hyperbolic discounting framework compared to their exponential discounting counterparts. This reveals the impatience of the decision makers. 

\subsection{Impact of present-bias }
\label{subsect3.3}
As has been revealed earlier that 
 under the stochastic quasi-hyperbolic discounting, the threshold for dividend payment is lower, which leads to earlier dividend payments, and earlier costly capital injections compared to the optimal scenario under the exponential discounting. To quantitatively assess the impact, we will compare  the cumulative amount of dividend payments, the cumulative amount of capital injections, and the value lost due to the adoption of the equilibrium strategy taken by the present biased decision maker under the long run  exponential interest rate.

 Recall that due to subjective present-biasedness, the decision maker now adopts the strategy $\pi^{0,b^*}$ instead of the optimal strategy $\pi^{0,b^E}$ without present biasedness.
 The value  under the equilibrium strategy is 
\begin{align}
\mathcal{P}^E(x;\pi^{0,b^*})
=&\mathrm{E}_{x}\bigg[\int^{\infty}_{0} e^{-\delta t}(d L^{0,b^*}{(t)}-\phi dR^{0,b^*}{(t)})\bigg]\nonumber\\
 =&V^E_{0,b^*}(x)\nonumber\\
 =&
\left\{
\small
\begin{aligned}
&
-
\overline{Z}_{\delta}(b^*-x)
-\frac{\psi^{\prime}(0+)}{\delta}
+\frac{Z_{\delta}(b^*-x)}{\delta W_{\delta}(b^*)}
\Big[Z_{\delta}(b^*)-\phi
\Big],& x\in[0,b^*],\\
&x-b^*+V^E_{0,b^*}(b),& x\in(b^*,\infty),\\
&\phi x+V^E_{0,b^*}(0), &x\in(-\infty,0),
\end{aligned}
\right.
\label{PLb-00}
\end{align}
which is the ``consistent" value of the company when adopting the equilibrium strategy. 
 Noting that the value without present-bias is $V^E(x)=\sup_{\pi}\mathcal{P}^E(x;\pi)$, which is greater than  $\mathcal{P}^E(x;\pi^{0,b^*})$,  we thus conclude that the behaviour issue has led to reduced value. The value lost can be quantified as
 $V^E_{0,b^E}(x)-V^E_{0,b^*}(x)$. We will provide numerical examples in the later sections and evaluate the lost value numerically. 
 
\section{``Optimality" of bailout/capital injections}\label{contionforcapitalinjections}

In the preceding sections, we addressed the optimal control problem within the context of a bailout scenario where forced capital injection prevents bankruptcy. What's more, the strategy $\pi^{0,b^\ast}=(R^0,L^{b^\ast})$ turned out to be an equilibrium strategy when capital injections are compulsory.
However, there might not always be sufficient economic incentive for a company to raise the costly capitals whenever needed. This section explores whether and when companies have sufficient economic incentives to raise capitals for business continuity. In scenarios without mandatory bailouts, the risk of bankruptcy arises when the surplus turns negative.
A natural question is: if there is no enforcement on bail-out/capital injections, is the strategy $\pi^{0,b^\ast}$ still an equilibrium dividend and capital injection strategy?  

We have discussed earlier that due to the transaction costs of raising capital, even if there are mechanisms to enforce compulsory capital injections to prevent bankruptcy, it is better not to inject capitals until the reserve reaches $0$.  
Let us look at the value associated with the equilibrium strategy $\pi^{0,b^\ast}$ when the reserve is $0$. 
It follows by Lemma \ref{Vb} that
\begin{align}
V_{0,b^\ast}(0)
=&-
\overline{Z}_{\delta+\lambda}(b^\ast)
-\frac{\psi^{\prime}(0+)}{\delta+\lambda}
+
\frac{\lambda\beta}{\delta+\lambda}\bigg[V_{0,b^\ast}^E(0)+
\int_{0}^{b^\ast}V_{0,b^\ast}^{E\,\prime}(y)
Z_{\delta+\lambda}(y){d}y\bigg].
\label{101023-2}
\end{align}
Intuitively, when $V_{0,b^*}(0)<0$, the strategy $\pi^{0,b^*}$ is not better than the strategy with no capital injections, in which case the value is always non-negative as it will contain the expected present of dividends only and dividends are always non-negative. Hence, it seems that a necessary conditions such that the strategy $\pi^{0,b^\ast}=(R^0,L^{b^\ast})$ is an equilibrium strategy can be $V_{0,b^{\ast}}(0)\ge 0$. We now want to explore whether this condition is a necessary and sufficient one for the strategy $\pi^{0,b^*}$ to be an equilibrium dividend and capital injection strategy when capital injections are not compulsory. 

Assume 
$V_{0,b^{\ast}}(0)\ge 0$. 
Let us define $\overline{\Pi}$ to represent the set of all the  strategies $\pi=(R,L)$ such that $L=\{L(t);t\ge 0\}$ and $R=\{R(t);t\ge 0\}$
are positive increasing cadl\'ag processes, adapted to  $(\mathcal{F}_t)_{t\ge 0}$,  the filtration generated by $\{X(t);t\ge 0\}$. Note that the set of admissible strategies we used previously, denoted by $\Pi$, is a subset of the expanded set $\overline{\Pi}$. This is because for any $\pi \in \Pi$, there is a constraint $X^\pi(t) \ge 0$ for all $t \ge 0$. The newly defined set $\overline{\Pi}$ has lifted the constraint on the positivity of the controlled process, allowing the company to inject no capital even at the risk of bankruptcy. The question we ask is  in the case $V_{0,b^{\ast}}(0)\ge 0$, whether the strategy $\pi^{0,b^*}$ is an equilibrium solution in  the larger set, $\overline{\Pi}$. That is, in the case $V_{0,b^{\ast}}(0)\ge 0$, does  $\mathcal{P}(x;\pi^{0,b^{\ast}},\pi^{0,b^{\ast}})=\sup_{\pi\in \overline{\Pi}}\mathcal{P}(x,\pi,\pi^{0,b^{\ast}})$ hold?

Define the bankruptcy/ruin time to be as
        \begin{align}
        \label{4.2}
T^\pi=\inf\{t\ge 0: U^\pi(t)< 0\}.
    \end{align}
By convention, we have $T^\pi=\infty$ if  $U^\pi(t)\ge 0$ for all $t\ge 0$. In particular, if $\pi\in\Pi$, then $T^\pi=\infty$. This is because under compulsory capital injections, which is the case for any strategy in $\Pi$, the controlled stochastic process will always stay at or above $0$ and hence, ruin never happens.  {\color{black}Recall that the natural filtration
$(\mathcal{F}_t)_{t\geq 0}$ generated by the spectrally positive L\'evy process $\{X(t),{t\geq 0}\}$ satisfies the usual conditions of right-continuity and completeness. Hence, $T^\pi$ defined by \eqref{4.2} is a stopping time.} Noting that there will be no cashflows after bankruptcy, one can write the  performance function associated with the strategy $\pi$ as
\begin{align}
&\mathcal{P}(x;\pi,\pi^{0,b^\ast})\nonumber\\
=&\mathrm{E}_{x}\bigg[\int^{\tau\wedge T^\pi}_{0} e^{-\delta t}(d L{(t)}-\phi dR{(t)})
+\beta e^{-\delta \tau} \mathcal{P}^E(U^{\pi}{(\tau)}; \pi^{0,b^\ast})\mathbf{1}_{\{  \tau< T^\pi\}}\bigg].
\label{4.3.q}
\end{align}
For any strategy $\pi=(R, L)$, construct a new strategy, denoted by $\hat{\pi}$, that applies the strategy $\pi$ during the time interval $[0,\tau\wedge T^\pi]$ and then employs the strategy $\pi^{0,b^\ast}$ during the time interval $(\tau\wedge T^\pi,\infty)$. Then, it holds that
\begin{align}
&\mathcal{P}(x;\hat{\pi},\pi^{0,b^\ast})\nonumber\\
=&\mathrm{E}_{x}\bigg[\int^{\tau\wedge T^\pi}_{0} e^{-\delta t}(d L{(t)}-\phi dR{(t)})
+\beta e^{-\delta \tau} \mathcal{P}^E(U^{\pi}{(\tau)}; \pi^{0,b^\ast})\mathbf{1}_{\{  \tau< T^\pi\}}
\nonumber\\
&\quad\quad +e^{-\delta T^\pi} \mathcal{P}(U^{\pi}(T^\pi);\pi^{0,b^\ast}, \pi^{0,b^\ast})\mathbf{1}_{\{  \tau> T^\pi\}}
\bigg]
\nonumber\\
=&\mathrm{E}_{x}\bigg[\int^{\tau\wedge T^\pi}_{0} e^{-\delta t}(d L{(t)}-\phi dR{(t)})
+\beta e^{-\delta \tau} \mathcal{P}^E(U^{\pi}{(\tau)}; \pi^{0,b^\ast})\mathbf{1}_{\{  \tau< T^\pi\}}
\nonumber\\
&\quad\quad +e^{-\delta T^\pi} V_{0,b^*}(0)\mathbf{1}_{\{  \tau> T^\pi\}}
\bigg],
\label{4.4.q}
\end{align}
where, in the second equality, we have used the facts that $U^{\pi}(T^\pi)=0$ and $\mathcal{P}(0;\pi^{0,b^\ast}, \pi^{0,b^\ast})=V_{0,b^*}(0)$.
Combining \eqref{4.3.q} and \eqref{4.4.q} and then using the fact that $V_{0,b^*}(0)\geq 0$, one can get
\begin{align}
\label{4.5.q}
    \mathcal{P}(x;\pi,\pi^{0,b^\ast})\leq \mathcal{P}(x;\hat{\pi},\pi^{0,b^\ast}), \quad \pi\in\overline{\Pi}.
\end{align}
In addition, by the construction of the strategy $\hat{\pi}$, one knows that $\hat{\pi}\in\Pi$. Therefore, by Theorem  \ref{16923-1}, we can obtain
\begin{align}
\label{4.6.q}
\mathcal{P}(x;\hat{\pi},\pi^{0,b^\ast})\leq
\mathcal{P}(x;\pi^{0,b^\ast},\pi^{0,b^\ast}).
\end{align}
Piecing together \eqref{4.5.q} and \eqref{4.6.q}, one arrives at
\begin{align}
\mathcal{P}(x;\pi,\pi^{0,b^\ast})\leq
\mathcal{P}(x;\pi^{0,b^\ast},\pi^{0,b^\ast}), \quad \pi\in\overline{\Pi}.
\end{align}
The arbitrariness of $\pi$ implies that \begin{align}
\sup_{\pi\in\overline{\Pi}}\mathcal{P}(x;{\pi},\pi^{0,b^\ast})\le \mathcal{P}(x;\pi^{0,b^\ast},\pi^{0,b^\ast}),\nonumber\end{align}
which, together with the fact that $\pi^{0,b^\ast}\in \Pi\subsetneq \overline{\Pi}$, yields
\begin{align}
\sup_{\pi\in\overline{\Pi}}
\mathcal{P}(x;{\pi},\pi^{0,b^\ast})= \mathcal{P}(x;\pi^{0,b^\ast},\pi^{0,b^\ast}).\nonumber\end{align}
Therefore, we can deduce that if $V_{0,b^{\ast}}(0)\ge 0$, the strategy $\pi^{0,b^\ast}$ remains an equilibrium strategy within the expanded admissible set $\overline{\Pi}$, which does not necessitate mandatory capital injections. This indicates that when $V_{0,b^{\ast}}(0)\ge 0$, it is advantageous to inject capital. In both scenarios-with or without compulsory bail-out requirements, the equilibrium solution is consistently to inject capital injections when necessary. 

On the contrary, if $V_{0,b^{\ast}}(0)< 0$, then the strategy $\pi^{0,b^\ast}$ cannot be an equilibrium strategy within the expanded admissible set $\overline{\Pi}$. Let $\pi^{b^\ast}=(R\equiv0,L^{b^\ast})$ be the dividend and capital injection strategy such that no capital injection is allowed and dividends are paid according to the barrier dividend strategy with barrier $b^\ast$. Then, it is known that the bankruptcy/ruin time $T^{\pi^{b^\ast}}$ (defined by \eqref{4.2} with $\pi$ replaced as $\pi^{b^\ast}$) is finite with probability 1. Similar to \eqref{4.3.q} and \eqref{4.4.q}, one has
\begin{align}
\mathcal{P}(x;\pi^{b^\ast},\pi^{0,b^\ast})&=\mathrm{E}_{x}\bigg[\int^{\tau\wedge T^{\pi^{b^\ast}}}_{0} e^{-\delta t}d L^{b^\ast}{(t)}
+\beta e^{-\delta \tau} \mathcal{P}^E(U^{\pi^{b^\ast}}{(\tau)}; \pi^{0,b^\ast})\mathbf{1}_{\{  \tau< T^{\pi^{b^\ast}}\}}\bigg],
\label{4.3.q.new}
\\
\mathcal{P}(x;\pi^{0,b^\ast},\pi^{0,b^\ast})
&=\mathrm{E}_{x}\bigg[\int^{\tau\wedge T^{\pi^{b^\ast}}}_{0} e^{-\delta t}d L^{b^\ast}{(t)}
+\beta e^{-\delta \tau} \mathcal{P}^E(U^{\pi^{b^\ast}}{(\tau)}; \pi^{0,b^\ast})\mathbf{1}_{\{  \tau< T^{\pi^{b^\ast}}\}}
\nonumber\\
&\quad\quad\,\,\,\,\,+e^{-\delta T^{\pi^{b^\ast}}} V_{0,b^*}(0)\mathbf{1}_{\{  \tau> T^{\pi^{b^\ast}}\}}
\bigg],
\label{4.4.q.new}
\end{align}
where we have used the fact that $U^{\pi^{b^\ast}}(t)=U^{\pi^{0,b^\ast}}(t)$ and $R^{0}(t)=0$ for all $t\in [0,T^{\pi^{b^\ast}}]$. Using \eqref{4.3.q.new}, \eqref{4.4.q.new}, and the fact that $V_{0,b^{\ast}}(0)< 0$ and $\mathrm{P}_{x}(T^{\pi^{b^\ast}}<\infty)=1$, one can deduce that
\begin{align}
\sup_{\pi\in\overline{\Pi}}\mathcal{P}(x;\pi,\pi^{0,b^\ast})
\geq 
\mathcal{P}(x;\pi^{b^\ast},\pi^{0,b^\ast})
>\mathcal{P}(x;\pi^{0,b^\ast},\pi^{0,b^\ast}).
\end{align}
Hence, $\pi^{0,b^\ast}$ is not an equilibrium strategy within the expanded admissible set $\overline{\Pi}$. 

In summary, the dividend and capital injection strategy $\pi^{0,b^\ast}$ remains an equilibrium strategy when bailout/capital injection is not compulsory, if and only if \begin{align}
  &  \overline{Z}_{\delta+\lambda}(b^\ast)
-\frac{\psi^{\prime}(0+)}{\delta+\lambda}
+
\frac{\lambda\beta}{\delta+\lambda}\bigg[V_{0,b^\ast}^E(0)+
\int_{0}^{b^\ast}V_{0,b^\ast}^{E\,\prime}(y)
Z_{\delta+\lambda}(y){d}y\bigg]
\ge 0, \label{111023-6}\end{align}
where $V_{0,b^\ast}^{E}(y)=
-
\overline{Z}_{\delta}(b^\ast-y)
-\frac{\psi^{\prime}(0+)}{\delta}
+\frac{Z_{\delta}(b^\ast-y)}{\delta W_{\delta}(b^\ast)}
\Big[Z_{\delta}(b^\ast)-\phi
\Big]$ for any $y\in[0,b^*]$.

\section{The Brownian Motion Case}\label{examples}
In this section, we will apply the  results obtained above to a classical case: a Brownian motion model.
Assume
$X{(t)}=x+\mu t+\sigma B{(t)}$, where $\mu$ and $\sigma$ ($\sigma>0$) are constants, and $\{B{(t)}, t\geq0\}$ is a standard Brownian motion. The Laplace exponent of $X$ is given by
\begin{eqnarray}
\psi(\theta):=\log\mathrm{E}(e^{\theta X(1)})=\frac{\sigma^2}{2}\theta^2+\mu\theta,\quad \theta\in (-\infty,\infty).
\nonumber
\end{eqnarray} Recall that for any $q\ge 0$, the $q$-scale functions $W_q(\cdot)$  is the unique strictly increasing and continuous function satisfying 
$$\int_{0}^{\infty}e^{-\theta y}W_q(y){d}y=\frac{1}{\psi(\theta)-q},\quad  \theta>\Phi_q,$$
where $\Phi_q$ is the largest solution of the equation $\psi(\theta)=q$, and $W_q(y)=0$ for $y<0$. We can determine the scale function: 
\begin{eqnarray}\label{W_q.example.1}
W_q(x)=\frac{1}{\sqrt{\mu^2+2\sigma^2q}}\left(e^{\frac{-\mu+\sqrt{\mu^2+2\sigma^2q}}{\sigma^2}x}-e^{\frac{-\mu-\sqrt{\mu^2+2\sigma^2q}}{\sigma^2}x}\right),\quad x\in[0,\infty).
\end{eqnarray}
Then the function $Z_q(x)$ defined in \eqref{8823-3} has the following representation: 
\begin{eqnarray}\label{Z_q.example.1}
Z_q(x)
\hspace{-0.3cm}&=&\hspace{-0.3cm}
\frac{q\sigma^2e^{\frac{-\mu+\sqrt{\mu^2+2\sigma^2q}}{\sigma^2}x}}{(-\mu+\sqrt{\mu^2+2\sigma^2q})\sqrt{\mu^2+2\sigma^2q}}+\frac{q\sigma^2e^{\frac{-\mu-\sqrt{\mu^2+2\sigma^2q}}{\sigma^2}x}}{(\mu+\sqrt{\mu^2+2\sigma^2q})\sqrt{\mu^2+2\sigma^2q}},\quad x\in[0,\infty).
\end{eqnarray}
Thus
\begin{eqnarray}
Z_q^\prime(x)
\hspace{-0.3cm}&=&\hspace{-0.3cm}
\frac{q\sigma^2\frac{-\mu+\sqrt{\mu^2+2\sigma^2q}}{\sigma^2}e^{\frac{-\mu+\sqrt{\mu^2+2\sigma^2q}}{\sigma^2}x}}{(-\mu+\sqrt{\mu^2+2\sigma^2q})\sqrt{\mu^2+2\sigma^2q}}+\frac{q\sigma^2\frac{-\mu-\sqrt{\mu^2+2\sigma^2q}}{\sigma^2}e^{\frac{-\mu-\sqrt{\mu^2+2\sigma^2q}}{\sigma^2}x}}{(\mu+\sqrt{\mu^2+2\sigma^2q})\sqrt{\mu^2+2\sigma^2q}}
\nonumber\\
\hspace{-0.3cm}&=&\hspace{-0.3cm}
\frac{q}{\sqrt{\mu^2+2\sigma^2q}}\left(e^{\frac{-\mu+\sqrt{\mu^2+2\sigma^2q}}{\sigma^2}x}-e^{\frac{-\mu-\sqrt{\mu^2+2\sigma^2q}}{\sigma^2}x}\right)
,\quad x\in[0,\infty).
\end{eqnarray}
Moreover, from \eqref{V.0bE.x} it follows that
\begin{align}
V_{0,b}^{E\,\prime}(x)=&
{Z}_{\delta}(b-x)
-\frac{Z_{\delta}^\prime(b-x)}{\delta W_{\delta}(b)}
\Big[Z_{\delta}(b)-\phi
\Big]\nonumber\\
=&\left(\frac{\delta\sigma^2}{(-\mu+\sqrt{\mu^2+2\sigma^2\delta})\sqrt{\mu^2+2\sigma^2\delta}}-\frac{Z_{\delta}(b)-\phi}{ W_{\delta}(b)\sqrt{\mu^2+2\sigma^2\delta}}\right)e^{\frac{-\mu+\sqrt{\mu^2+2\sigma^2\delta}}{\sigma^2}(b-x)}\nonumber\\
&+\left(\frac{\delta\sigma^2}{(\mu+\sqrt{\mu^2+2\sigma^2\delta})\sqrt{\mu^2+2\sigma^2\delta}}+\frac{Z_{\delta}(b)-\phi}{ W_{\delta}(b)\sqrt{\mu^2+2\sigma^2\delta}}\right)e^{\frac{-\mu-\sqrt{\mu^2+2\sigma^2\delta}}{\sigma^2}(b-x)},\quad x\in[0,b].\label{mono.cond.2}
\end{align}
Write $q:=\delta+\lambda$.
Recall that $\ell(x)=Z_{q}(x)
-\lambda\beta\int_{0}^{x}V_{0,x}^{E\,\prime}(y)W_{q}(y){d}y-\phi.$  By noticing that $V_{0,b}^{E\,\prime}(\cdot)$ and  $W_q(\cdot)$ are both sum of exponential functions, we can easily find an explicit expression for $\ell(x)$, which is also a sum of exponential function of $x$:
\begin{eqnarray}
\ell(x)
\hspace{-0.3cm}&=&\hspace{-0.3cm}
\frac{q\sigma^2e^{\frac{-\mu+\sqrt{\mu^2+2\sigma^2q}}{\sigma^2}x}}{(-\mu+\sqrt{\mu^2+2\sigma^2q})\sqrt{\mu^2+2\sigma^2q}}+\frac{q\sigma^2e^{\frac{-\mu-\sqrt{\mu^2+2\sigma^2q}}{\sigma^2}x}}{(\mu+\sqrt{\mu^2+2\sigma^2q})\sqrt{\mu^2+2\sigma^2q}}
\nonumber\\
\hspace{-0.3cm}&&\hspace{-0.3cm}
-\lambda\beta\Big(\frac{\delta\sigma^2}{(-\mu+\sqrt{\mu^2+2\sigma^2\delta})\sqrt{\mu^2+2\sigma^2\delta}}
-\frac{Z_{\delta}(x)-\phi}{W_{\delta}(x)\sqrt{\mu^2+2\sigma^2\delta}}\Big)
\nonumber\\
\hspace{-0.3cm}&&\hspace{-0.3cm}
\times\frac{1}{\sqrt{\mu^2+2\sigma^2q}}\frac{\sigma^2}{-\sqrt{\mu^2+2\sigma^2\delta}+\sqrt{\mu^2+2\sigma^2q}}\left(e^{\frac{-\mu+\sqrt{\mu^2+2\sigma^2q}}{\sigma^2}x}-e^{\frac{-\mu+\sqrt{\mu^2+2\sigma^2\delta}}{\sigma^2}x}\right)
\nonumber\\
\hspace{-0.3cm}&&\hspace{-0.3cm}
-\lambda\beta\Big(\frac{\delta\sigma^2}{(\mu+\sqrt{\mu^2+2\sigma^2\delta})\sqrt{\mu^2+2\sigma^2\delta}}
+\frac{Z_{\delta}(x)-\phi}{W_{\delta}(x)\sqrt{\mu^2+2\sigma^2\delta}}\Big)
\nonumber\\
\hspace{-0.3cm}&&\hspace{-0.3cm}
\times\frac{1}{\sqrt{\mu^2+2\sigma^2q}}\frac{\sigma^2}{\sqrt{\mu^2+2\sigma^2\delta}+\sqrt{\mu^2+2\sigma^2q}}\left(e^{\frac{-\mu+\sqrt{\mu^2+2\sigma^2q}}{\sigma^2}x}-e^{\frac{-\mu-\sqrt{\mu^2+2\sigma^2\delta}}{\sigma^2}x}\right)
\nonumber\\
\hspace{-0.3cm}&&\hspace{-0.3cm}
-\lambda\beta\Big(\frac{\delta\sigma^2}{(-\mu+\sqrt{\mu^2+2\sigma^2\delta})\sqrt{\mu^2+2\sigma^2\delta}}
-\frac{Z_{\delta}(x)-\phi}{W_{\delta}(x)\sqrt{\mu^2+2\sigma^2\delta}}\Big)
\nonumber\\
\hspace{-0.3cm}&&\hspace{-0.3cm}
\times\frac{1}{\sqrt{\mu^2+2\sigma^2q}}\frac{\sigma^2}{\sqrt{\mu^2+2\sigma^2\delta}+\sqrt{\mu^2+2\sigma^2q}}\left(e^{\frac{-\mu-\sqrt{\mu^2+2\sigma^2q}}{\sigma^2}x}-e^{\frac{-\mu+\sqrt{\mu^2+2\sigma^2\delta}}{\sigma^2}x}\right)
\nonumber\\
\hspace{-0.3cm}&&\hspace{-0.3cm}
-\lambda\beta\Big(\frac{\sigma^2}{(\mu+\sqrt{\mu^2+2\sigma^2\delta})\sqrt{\mu^2+2\sigma^2\delta}}
+\frac{Z_{\delta}(x)-\phi}{W_{\delta}(x)\sqrt{\mu^2+2\sigma^2\delta}}\Big)
\nonumber\\
\hspace{-0.3cm}&&\hspace{-0.3cm}
\times\frac{1}{\sqrt{\mu^2+2\sigma^2q}}\frac{\sigma^2}{\sqrt{\mu^2+2\sigma^2\delta}-\sqrt{\mu^2+2\sigma^2q}}\left(e^{\frac{-\mu-\sqrt{\mu^2+2\sigma^2q}}{\sigma^2}x}-e^{\frac{-\mu-\sqrt{\mu^2+2\sigma^2\delta}}{\sigma^2}x}\right)-\phi.
\end{eqnarray}
By finding the first positive solution to $\ell(x)=0$, we can identify $b^\ast$.

\begin{figure}[h!tb]
\centering
\includegraphics[width=3.2in]{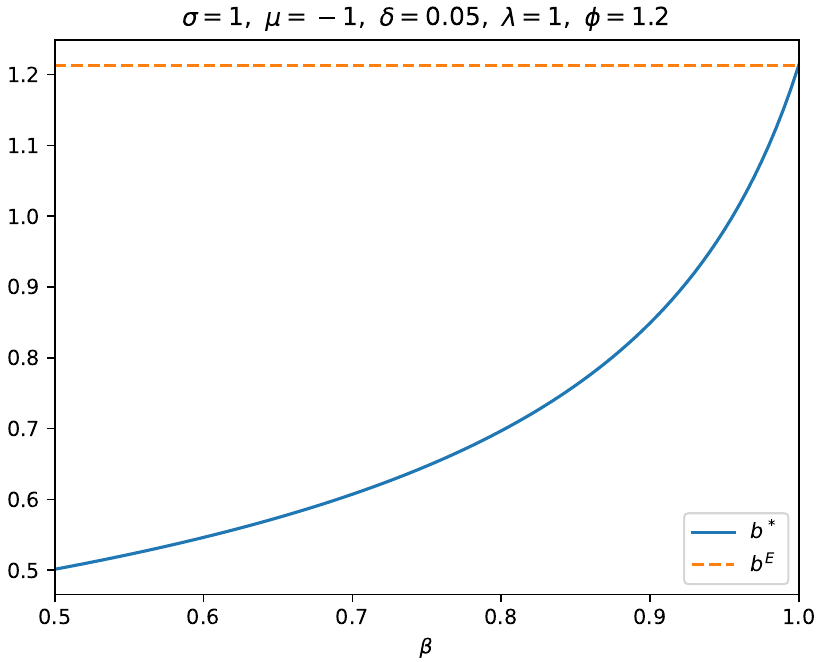}
\caption{Example 4.1: The optimal dividend barriers when $\beta$ varies}
\label{figure1}
\end{figure}

Now we proceed with numerical illustrations. First we set $\mu=-1$, $\sigma=1$, and $\delta=5\%$. In Figure \ref{figure1}, 
the barriers of the equilibrium dividend strategies are depicted  for various $\beta$ when $\lambda=1$ and $\phi=1.2$.  We can see that as $\beta$ increases, the  threshold $b^\ast$ (the solid curve) of the equilibrium dividend strategy increases. This is because as $\beta$ increases, the discounting for the future cashflows is getting closer to the exponential discounting, which  means the decision makers are more patient in the cases with a larger $\beta$, and when $\beta=1$, they are as patient as the counterpart in exponential discounting, and hence  identical dividend payment thresholds in the stochastic quasi-hyperbolic discounting case with $\beta=1$ and the exponential discounting scenario (represented by the dotted line).

\begin{figure}[h!tb]
\centering
\includegraphics[width=3.2in]{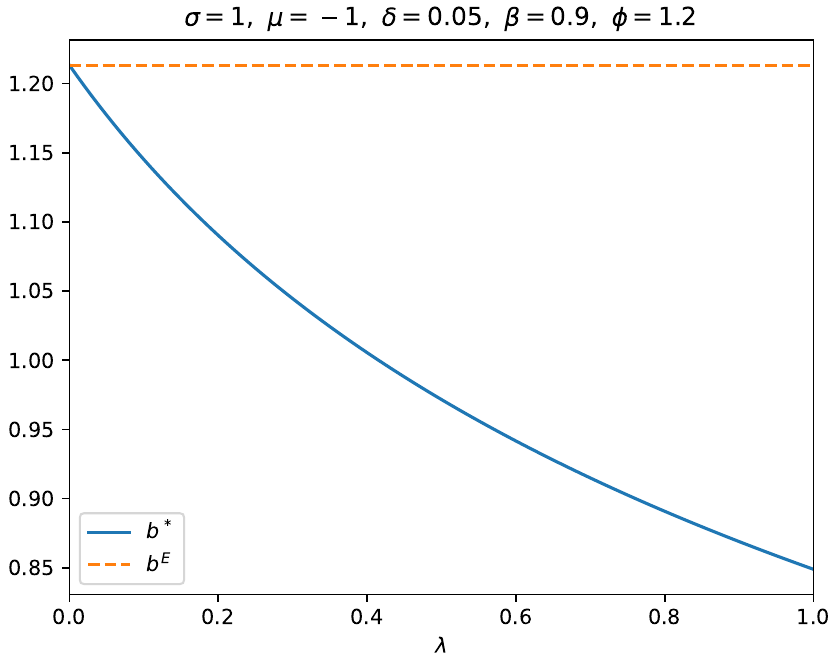}
\caption{Example 4.1: The optimal dividend barriers when $\lambda$ varies}
\label{figure2-1}
\end{figure}

Figure \ref{figure2-1} 
depicts the barriers of the equilibrium strategies for various $\lambda$ when $\beta=0.9$ and $\phi=1.2$. As observed, 
when the arrival intensity of the future periods increases, the dividend payment thresholds getting smaller. This is because the future arrives earlier with higher likelihood when $\lambda$ increases, which means that the decision makers are more impatient and desire for dividend payments earlier.

\begin{figure}[h!tb]
\centering
\includegraphics[width=3.2in]{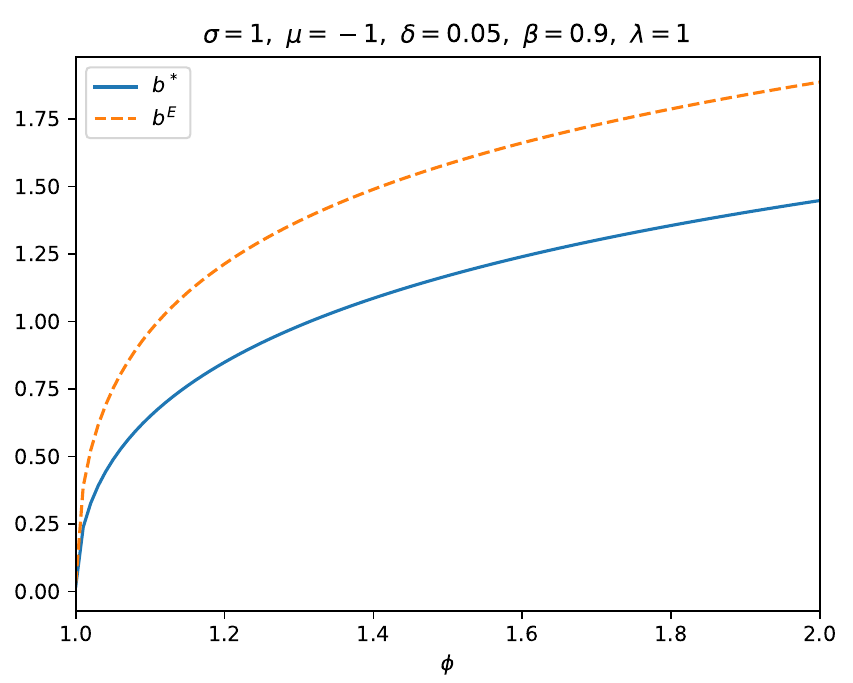}
\caption{Example 4.1: The optimal dividend barriers when $\phi$ varies}
\label{figure3-0}
\end{figure}
We plot the barriers of the equilibrium dividend strategies for various $\phi$ (the capital injection cost) in Figure \ref{figure3-0}.
In both the exponential and stochastic quasi-hyperbolic discounting cases, when capital injection cost is $0$ (i.e., $\phi=1$), the dividend barrier is $0$. This means that when capital injections are cost free,  it is always optimal to pay all the available surplus  out as dividends  and to inject capitals whenever necessarily to keep the business running. When the cost of capital injection increases in both scenarios, the company raises the barrier to reduce dividend distribution, thereby decreasing the frequency of needing to raise costly capital.

\begin{figure}[h!tb]
\centering
\includegraphics[width=3.2in]{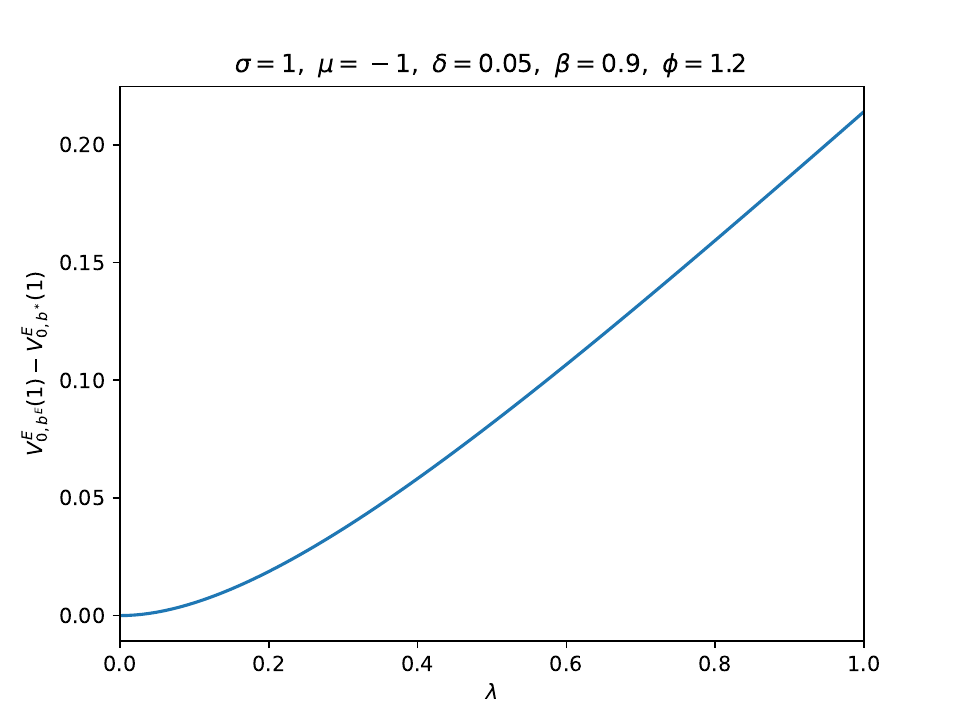}
\quad
\includegraphics[width=3.2in]{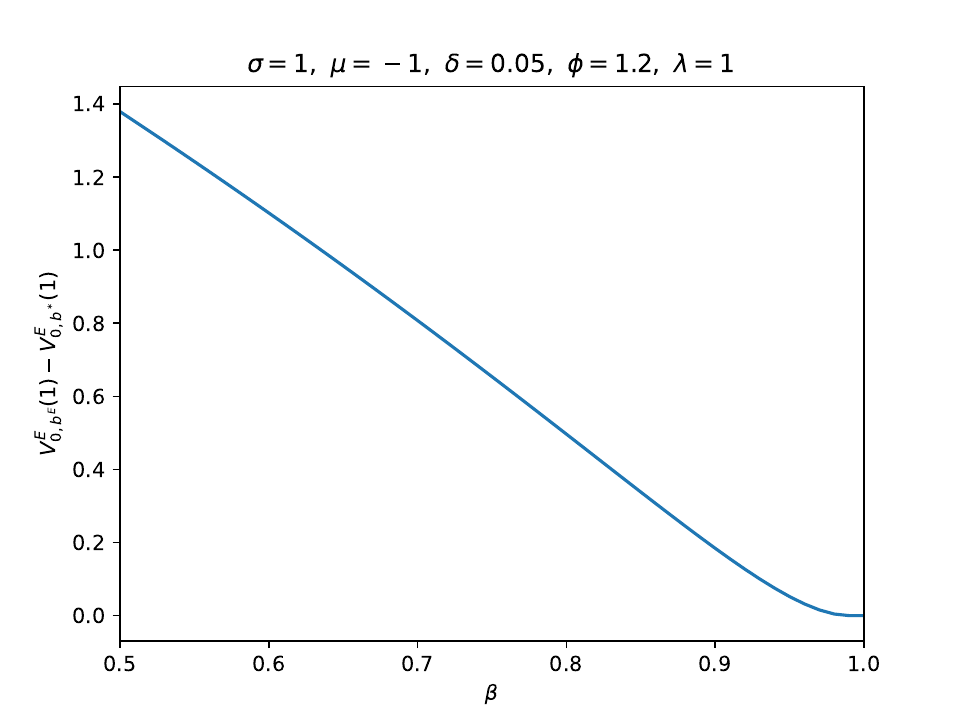}
\caption{Example 4.1: Losses}
\label{91223-1}
\end{figure}

To comprehend the impact of behavior issues on the incurred loss, we depict the loss (the disparity in values between cases with and without behavior issues) for the case with an initial surplus $1$ in Figure \ref{91223-1}. 
The observed trend reveals a consistent presence of loss attributable to behavior issues, with the magnitude increasing as the intensity parameter $\lambda$ rises or the additional discount factor $\beta$  on the future cashflows decreases. This observation implies that higher levels of present-bias lead to more pronounced losses.

%

\begin{figure}[htbp]
\label{figure4}
\centering
\includegraphics[width=4in]{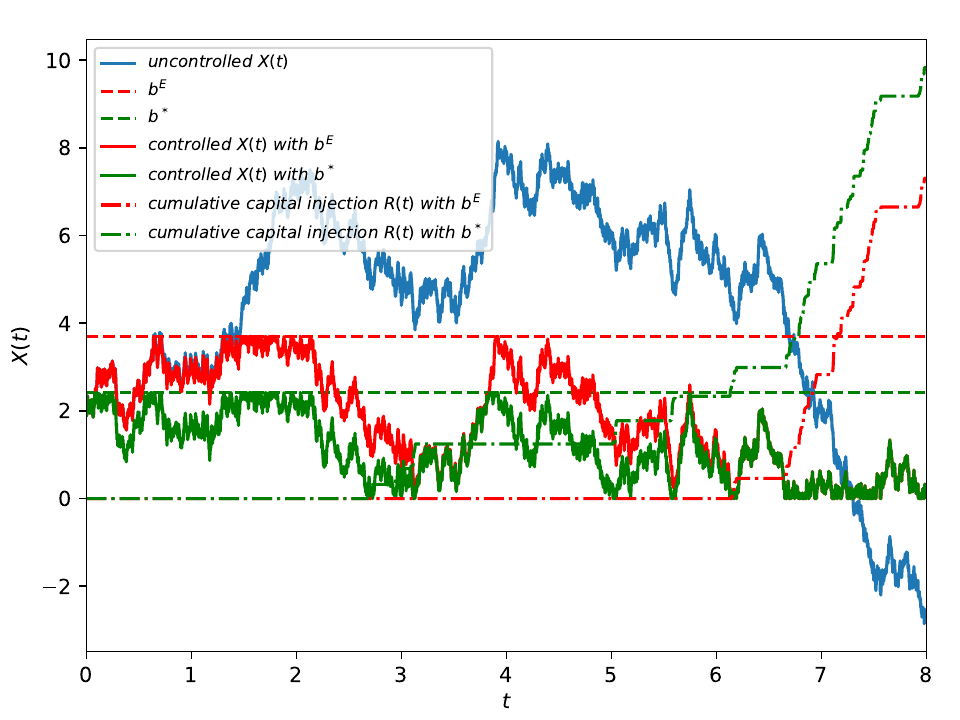}
\caption{Example 4.1: The uncontrolled process and the controlled process $X_t$ with capital injection barrier $0$ and dividend payment barrier $b^E$ and $b^*$ when $x=2,\,\mu=-0.2,\,\sigma=2,\,\delta=0.05,\,\beta=0.9,\,\lambda=1,\,\phi=1.2,$ and $V_{0,b^*}(0)<0$}

\end{figure}

We now simulate a path for the uncontrolled L\'evy process, and then the paths of the corresponding optimally controlled processes in the stochastic quasi-hyperbolic discounting case under the equilibrium strategy, and the exponential discounting case under the optimal strategy, respectively. In Figure \ref{figure4}, the blue solid curve represents a path of the uncontrolled L\'evy  process starting with an initial value $2$, the red solid curve depicts the corresponding path if it
is optimally controlled by the optimal strategy (which has an optimal capital injection barrier $0$, and an optimal barrier at the level, $b^E$, depicted by the red dotted line) in the exponential discounting case. The green curve represents the corresponding path of the controlled stochastic process under the equilibrium strategy, which is a double barrier strategy with capital injection barrier $0$ and dividend barrier, $b^*$ (depicted by the dotted green live), in the stochastic quasi-hyperbolic case.  We can see that in the stochastic quasi-hyperbolic case,  the payment barrier is much lower than the exponential discounting case and as a result, the controlled surplus is lower as well, which leads to earlier capital injections. We can observe that the cumulative capital injections in the stochastic quasi-hyperbolic case (the green case) is higher than those in the exponential case (the red case).  \\

\begin{figure}[H]
\centering
\includegraphics[width=3.2in]{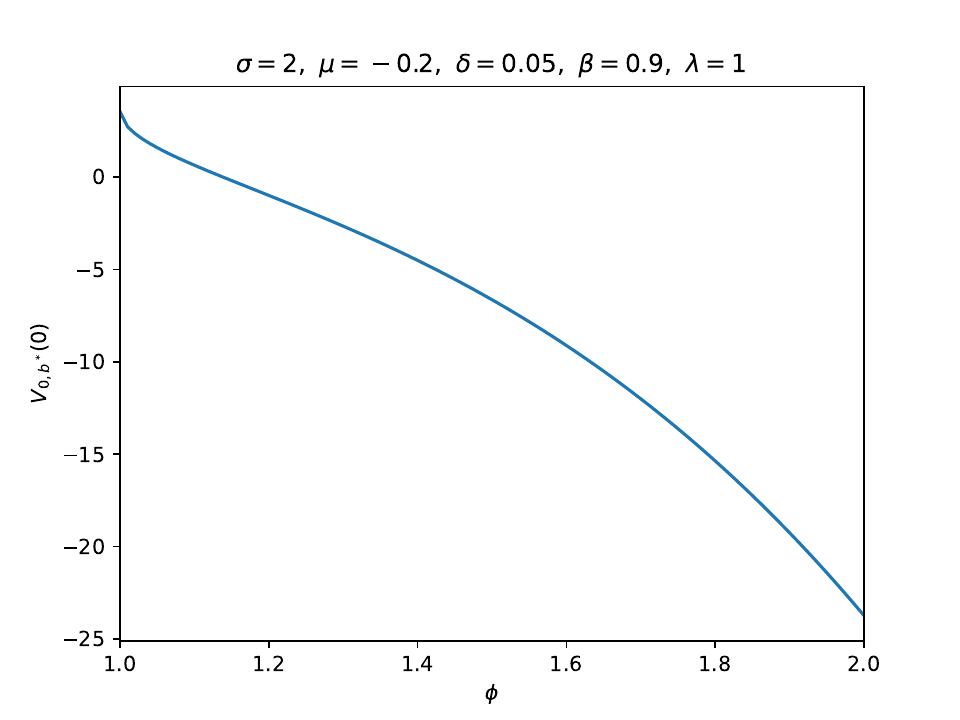}
\caption{Example 4.1: $V_{0,b^*}(0)$ vs $\phi$ \label{fig15124-1}}

\end{figure}

We have shown in Section \ref{contionforcapitalinjections} that the condition for the optimality of injecting capitals is $V_{0,b^*}(0)\ge 0$.
A graph depicting $V_{0,b^*}(0)$ in relation to $\phi$ (refer to Figure \ref{fig15124-1}) reveals that when the capital injection additional cost factor $\phi$ is small and close to $1$, it becomes optimal to inject capital. Notably, $\phi=1$ indicates no additional cost.\\

\section{The Jump-diffusion Case}
\label{jumpdiffusion}

In this section,  we consider the jump diffusion model with exponential jump size distributions. Jump diffusion models are extensively used in finance and actuarial science.  Let $X{(t)}=x+\mu t+\sigma B{(t)}+\sum_{i=1}^{N{(t)}}e_i$, where $x$ is the initial endowment, constant $\mu$ is the drift coefficient,  $\sigma>0$ is the dispersion parameter, $\{B{(t)},t\geq0\}$ is a standard Brownian motion,  $\{e_i;i\geq1\}$ are independent and identically distributed exponential random variables with mean $1/\eta$, and $\{N{(t)},t\geq0\}$ is an independent Poisson process with intensity $p$. Then the Laplace exponent of $X$ is given by
\begin{eqnarray}
\psi(\theta):=\log\mathrm{E}(e^{\theta X(1)})=\frac{\sigma^2}{2}\theta^2+\mu\theta+p\big[\frac{\eta}{\eta-\theta}-1\big],\quad \theta\in (-\infty,\eta).
\nonumber
\end{eqnarray}
Recall $q=\delta+\lambda$.
Further define $W_q(x)=0$ for $x<0$.
The scale function $W_q(x)$ reads as
\begin{eqnarray}\label{W_q.example.2}
W_q(x)=\frac{e^{\vartheta_q x}}{\psi^{\prime}(\vartheta_q)}+\frac{e^{\theta_qx}}{\psi^{\prime}(\theta_q)}+\frac{e^{\Phi_q x}}{\psi^{\prime}(\Phi_q)},\quad x\in [0,\infty),
\end{eqnarray}
where $\vartheta_q$, $\theta_q$ and $\Phi_q$ denote the three roots of $\psi(\theta)=q$ such that $\vartheta_{q} < 0 < \theta_{q} <\eta<\Phi_{q}$ and $\psi^{\prime}(\vartheta_q)<0$, $\psi^{\prime}(\theta_q)>0$ and $\psi^{\prime}(\Phi_q)>0$.
Then the scale function $Z_q(x)$ can be rewritten as
\begin{eqnarray}\label{Z_q.example.2}
Z_q(x)=1+q\int_0^xW_q(y){d}y=\frac{qe^{\vartheta_q x}}{\psi^{\prime}(\vartheta_q)\vartheta_q}+\frac{qe^{\theta_qx}}{\psi^{\prime}(\theta_q)\theta_q}+\frac{qe^{\Phi_q x}}{\psi^{\prime}(\Phi_q)\Phi_q},\quad x\in [0,\infty).
\end{eqnarray}
Combining \eqref{mono.cond.2}, \eqref{W_q.example.2} and \eqref{Z_q.example.2}, we obtain
\begin{eqnarray}
\ell(x)
\hspace{-0.3cm}&=&\hspace{-0.3cm}
\bigg[\frac{q}{\vartheta_q}-\Big(\frac{\phi-Z_{\delta}(x)}{W_{\delta}(x)}+\frac{\delta}{\vartheta_{\delta}}\Big)\frac{\lambda\beta}{\psi^{\prime}(\vartheta_{\delta})(\vartheta_q-\vartheta_{\delta})}-\Big(\frac{\phi-Z_{\delta}(x)}{W_{\delta}(x)}+\frac{\delta}{\theta_{\delta}}\Big)\frac{\lambda\beta}{\psi^{\prime}(\theta_{\delta})(\vartheta_q-\theta_{\delta})}
\nonumber\\
\hspace{-0.3cm}&&\hspace{-0.3cm}
-\Big(\frac{\phi-Z_{\delta}(x)}{W_{\delta}(x)}+\frac{\delta}{\Phi_{\delta}}\Big)\frac{\lambda\beta}{\psi^{\prime}(\Phi_{\delta})(\vartheta_q-\Phi_{\delta})}\bigg]\frac{e^{\vartheta_qx}}{\psi^{\prime}(\vartheta_q)}
\nonumber\\
\hspace{-0.3cm}&&\hspace{-0.3cm}
+\bigg[\frac{q}{\theta_q}-\Big(\frac{\phi-Z_{\delta}(x)}{W_{\delta}(x)}+\frac{\delta}{\vartheta_{\delta}}\Big)\frac{\lambda\beta}{\psi^{\prime}(\vartheta_{\delta})(\theta_q-\vartheta_{\delta})}-\Big(\frac{\phi-Z_{\delta}(x)}{W_{\delta}(x)}+\frac{\delta}{\theta_{\delta}}\Big)\frac{\lambda\beta}{\psi^{\prime}(\theta_{\delta})(\theta_q-\theta_{\delta})}
\nonumber\\
\hspace{-0.3cm}&&\hspace{-0.3cm}
-\Big(\frac{\phi-Z_{\delta}(x)}{W_{\delta}(x)}+\frac{\delta}{\Phi_{\delta}}\Big)\frac{\lambda\beta}{\psi^{\prime}(\Phi_{\delta})(\theta_q-\Phi_{\delta})}\bigg]\frac{e^{\theta_qx}}{\psi^{\prime}(\theta_q)}
\nonumber\\
\hspace{-0.3cm}&&\hspace{-0.3cm}
+\bigg[\frac{q}{\Phi_q}-\Big(\frac{\phi-Z_{\delta}(x)}{W_{\delta}(x)}+\frac{\delta}{\vartheta_{\delta}}\Big)\frac{\lambda\beta}{\psi^{\prime}(\vartheta_{\delta})(\Phi_q-\vartheta_{\delta})}-\Big(\frac{\phi-Z_{\delta}(x)}{W_{\delta}(x)}+\frac{\delta}{\theta_{\delta}}\Big)\frac{\lambda\beta}{\psi^{\prime}(\theta_{\delta})(\Phi_q-\theta_{\delta})}
\nonumber\\
\hspace{-0.3cm}&&\hspace{-0.3cm}
-\Big(\frac{\phi-Z_{\delta}(x)}{W_{\delta}(x)}+\frac{\delta}{\Phi_{\delta}}\Big)\frac{\lambda\beta}{\psi^{\prime}(\Phi_{\delta})(\Phi_q-\Phi_{\delta})}\bigg]\frac{e^{\Phi_qx}}{\psi^{\prime}(\Phi_q)}-\phi
\nonumber\\
\hspace{-0.3cm}&&\hspace{-0.3cm}
+\lambda\beta\Big(\frac{\phi-Z_{\delta}(x)}{W_{\delta}(x)\psi^{\prime}(\vartheta_{\delta})}+\frac{\delta}{\vartheta_{\delta}\psi^{\prime}(\vartheta_{\delta})}\Big)\Big[\frac{1}{\psi^{\prime}(\vartheta_q)(\vartheta_q-\vartheta_{\delta})}+\frac{1}{\psi^{\prime}(\theta_q)(\theta_q-\vartheta_{\delta})}+\frac{1}{\psi^{\prime}(\Phi_q)(\Phi_q-\vartheta_{\delta})}\Big]e^{\vartheta_{\delta}x}
\nonumber\\
\hspace{-0.3cm}&&\hspace{-0.3cm}
+\lambda\beta\Big(\frac{\phi-Z_{\delta}(x)}{W_{\delta}(x)\psi^{\prime}(\theta_{\delta})}+\frac{\delta}{\theta_{\delta}\psi^{\prime}(\theta_{\delta})}\Big)\Big[\frac{1}{\psi^{\prime}(\vartheta_q)(\vartheta_q-\theta_{\delta})}+\frac{1}{\psi^{\prime}(\theta_q)(\theta_q-\theta_{\delta})}+\frac{1}{\psi^{\prime}(\Phi_q)(\Phi_q-\theta_{\delta})}\Big]e^{\theta_{\delta}x}
\nonumber\\
\hspace{-0.3cm}&&\hspace{-0.3cm}
+\lambda\beta\Big(\frac{\phi-Z_{\delta}(x)}{W_{\delta}(x)\psi^{\prime}(\Phi_{\delta})}+\frac{\delta}{\Phi_{\delta}\psi^{\prime}(\Phi_{\delta})}\Big)\Big[\frac{1}{\psi^{\prime}(\vartheta_q)(\vartheta_q-\Phi_{\delta})}+\frac{1}{\psi^{\prime}(\theta_q)(\theta_q-\Phi_{\delta})}+\frac{1}{\psi^{\prime}(\Phi_q)(\Phi_q-\Phi_{\delta})}\Big]e^{\Phi_{\delta}x}.\nonumber
\end{eqnarray}
The dividend barrier of the equilibrium strategy can be determined by finding the positive root of $\ell(x)$.

We now present some numerical results by setting $\mu=-1$, $\sigma=2$, $\delta=5\%$, $\beta=0.9$ and $\phi=1.2$, and then varies some of the parameters one by one separately by fixing the other parameters to investigate the sensitivity of the $b^\ast$ to each of those parameters. 
\begin{figure}[h!tb]
\centering
\includegraphics[width=3.2in]{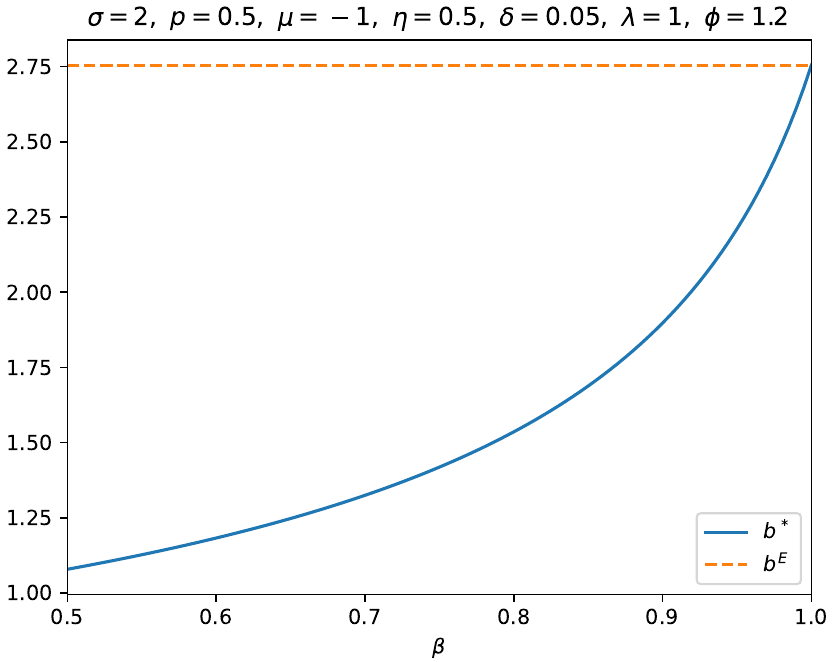}
\caption{Example 4.2: The optimal barriers when $\beta$ varies }
\label{figure7}
\end{figure}

\begin{figure}[h!tb]
\centering
\includegraphics[width=3.2in]{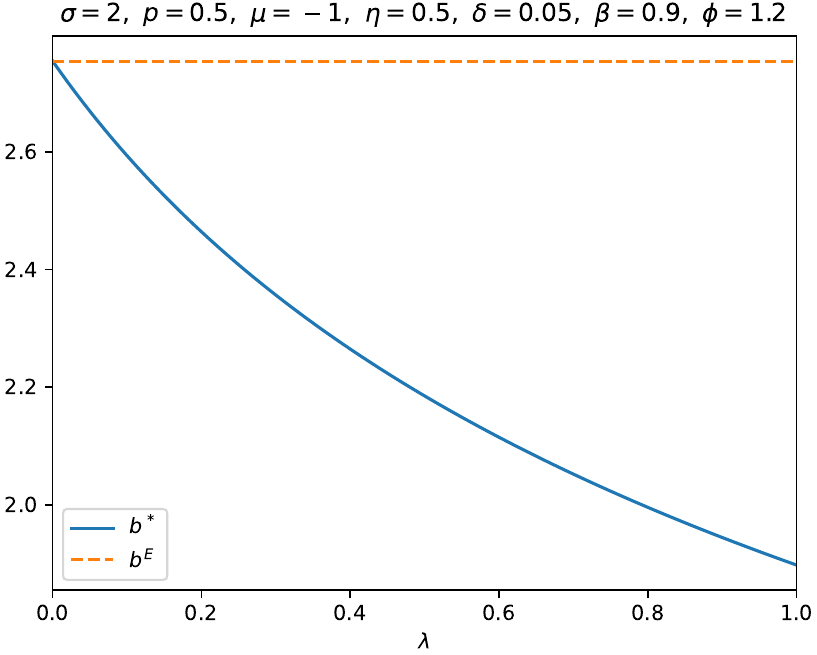}
\caption{Example 4.2: The optimal barriers when $\lambda$ varies }
\label{figure8}
\end{figure}

\begin{figure}[h!tb]
\centering
\includegraphics[width=3.4in]{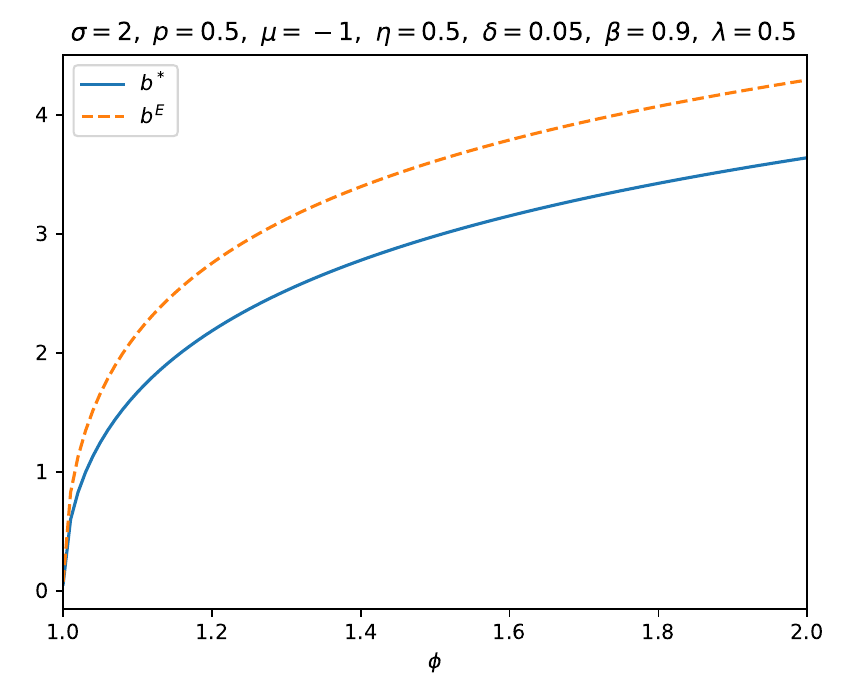}
\caption{Example 4.2: The optimal barriers when $\phi$ varies }
\label{figure9}
\end{figure}
Similar to the pure diffusion case, in the jump diffusion framework, the  threshold $b^\ast$ (the solid curves in Figures \ref{figure7} and \ref{figure8})
of the equilibrium dividend strategy  is lower for smaller $\beta$ and larger $\lambda$, respectively,  which implies that dividends are paid earlier and more dividends are paid out when there is higher impatience.
Regarding the impact of the cost of capital injections, in Figure \ref{figure9} 
we observe similar phenomenon as in the pure diffusion case. When $\phi=1$ (no cost): it is always optimal to pay all the available surplus  out as dividends  and then to inject capitals whenever necessarily, and the company lifts the dividend payment barriers when the cost of capital injection increases in order to reduce  capital injections.

\begin{figure}[h!tb]
\centering
\includegraphics[width=3.2in]{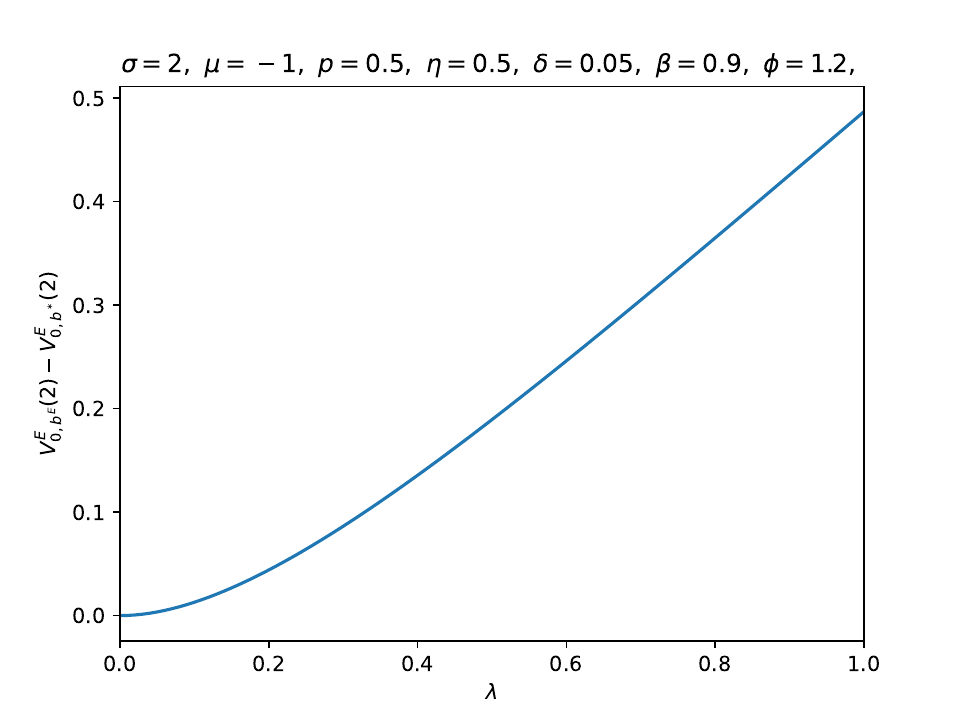}
\quad
\includegraphics[width=3.2in]{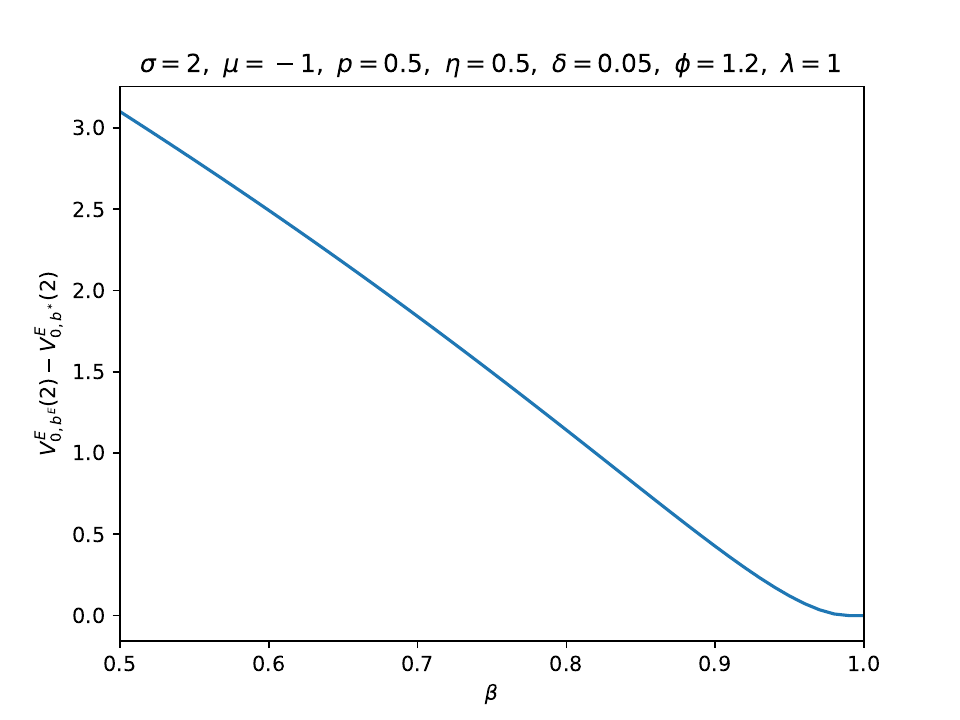}
\caption{Example 4.2: Losses}
\label{91223-2}
\end{figure}

The losses  resulting from present-bias are illustrated in Figure \ref{91223-2} (for the case with an  initial surplus $2$). 
Evidently, as the level of present-bias increases (reflected by higher values of $\lambda$ or lower values of $\beta$), the incurred losses also escalate.

When $0<\beta<1$, the optimal barrier, $b{^\ast}$, is always smaller than $b^E$, and $b{^\ast}$ becomes smaller as $\beta$ moves away from below $1$. As the strategy prescribes to pay dividends once the barrier is reached, a lower barrier implies that payments will start earlier and higher than the decision maker would like using exponential discounting. This will also leads to earlier capital injection and higher amount of cumulative injections, which is not desirable sometimes due to the transaction cost.

\begin{figure}[h!tb]
\centering
\includegraphics[width=4in]{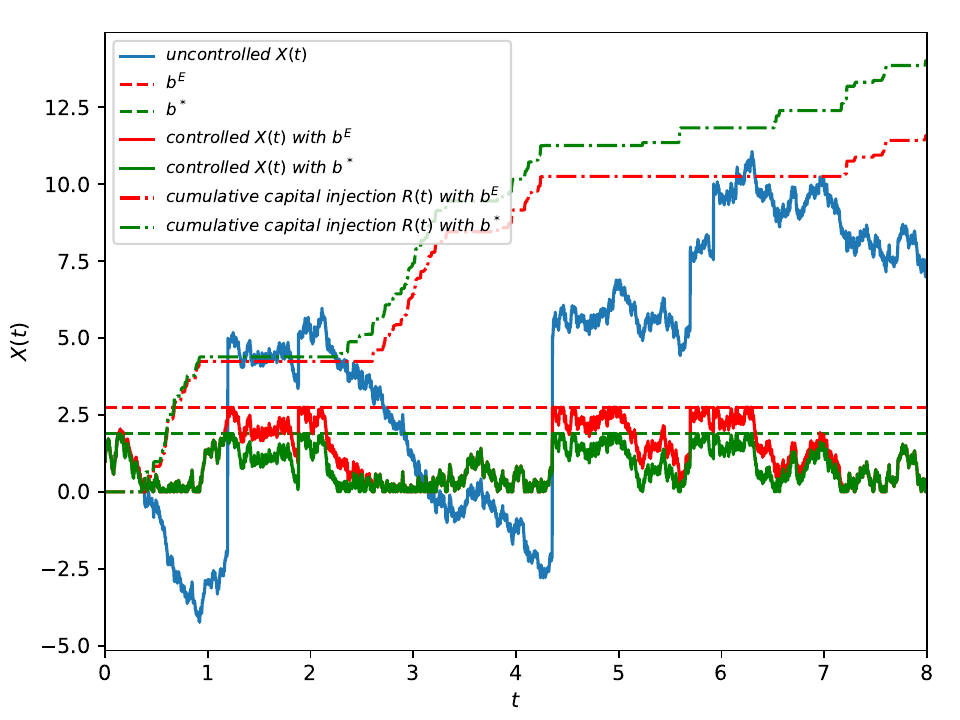}
\caption{Example 4.2: The uncontrolled process and the controlled process $X_t$ with the dividend payment threshold being $b^E$ and $b^*$ when $x=1,\,\mu=-1,\,p=0.5,\,\sigma=2,\,\delta=0.05,\,\beta=0.9,\,\lambda=1,\,\phi=1.2$,\text{ and }$V_{0,b^*}(0)<0$}
\label{figure10}
\end{figure}

Figure \ref{figure10} 
provides  simulated trajectories of a stochastic process  where it is never optimal to inject capitals at all.
We simulate a trajectory for the uncontrolled L\'evy process, followed by the trajectories of the optimally controlled processes in two scenarios: the stochastic quasi-hyperbolic discounting case under the equilibrium strategy, and the exponential discounting case under the optimal strategy. The blue solid curve represents the trajectory of the uncontrolled L\'evy process, initiated with an initial value of $1$. The red solid curve illustrates the corresponding trajectory when the stochastic process is optimally controlled in the exponential discounting case (the case with no present-bias), while the green curve represent the present-biased case. 
 We observe that there are earlier and more capital injections when the decision  makers are present-biased (in the stochastic quasi-hyperbolic discounting case): the green dotted curve is larger than the red one.  

\begin{figure}[h!tb]
\centering
\includegraphics[width=3.2in]{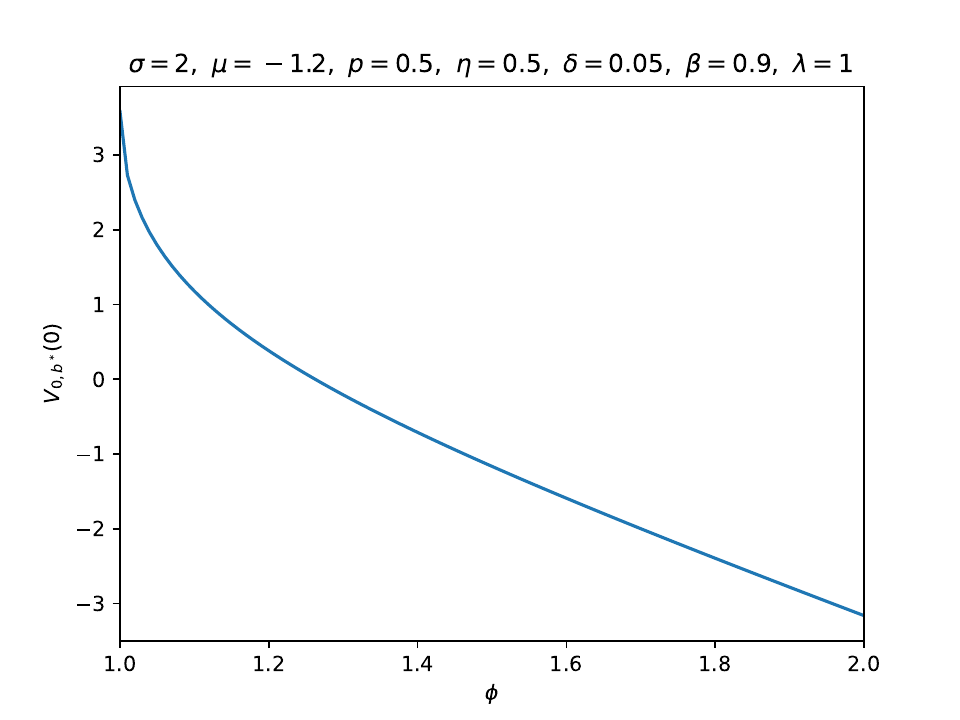}
\caption{Example 4.2: $V_{0,b^*}(0)$ vs $\phi$}
\label{fig15124-2}
\end{figure}

Recall from Section \ref{contionforcapitalinjections} we know that the condition for the optimality of injecting capitals is $V_{0,b^*}(0)\ge 0$.
A plot of  $V_{0,b^*}(0)$ versus $\phi$ (refer to Figure \ref{fig15124-2}) shows that when the capital injection additional cost factor $\phi$ is small and close to $1$ (smaller than $1.25$), it is optimal to inject capital. Note $\phi=1$ means no additional cost.\\

\section{Conclusion}\label{conclusion}
This study integrates stochastic quasi-hyperbolic discounting into the optimization of earning/retaining and capital injection control under spectrally positive L\'evy processes. This research models this as an intra-personal game, and derives a stationary Markov perfect equilibrium (MPE). Notably, the widely used double barrier strategies encompass an MPE. The outcomes highlight a tendency for impatient dividend payments, leading firms to initiate dividends earlier (reflected in the lower barrier) than exponential discounting scenarios. The consequence is reduced profit (value) and more costly capital injections.

The approach presented in this paper can readily be expanded to address diverse bailout challenges involving various L\'evy processes, such as spectrally negative ones. Additionally, with suitable adjustments, it can be applied to tackle different value maximization problems under L\'evy processes.

\section{Appendix}
\renewcommand{\theequation}{A-\arabic{equation}}
\setcounter{equation}{0}
\textbf{Proof for Lemma \ref{lem2.4}}. (i) It follows by  \eqref{V.0bE.x} that
\begin{eqnarray}
V_{0,b}^{E\,\prime}(x)
\hspace{-0.3cm}&=&\hspace{-0.3cm}
Z_{\delta}(b-x)+\frac{\phi-Z_{\delta}(b)}{W_{\delta}(b)}W_{\delta}(b-x),\quad
 x\in(0,b),\label{9823-1}
\end{eqnarray}
which, together with the strictly decreasing property of $Z_{\delta}(b-x)$ and $W_\delta(b-x)$ in $x$, and the fact that $\phi-Z_{\delta}(b)\geq 0$ for $b\le b^E$ (by the definition of $b^E$ in Definition \ref{bE.def}), implies that  ${V_{0,b}^{E\,\prime}}(x)$ is strictly decreasing on $(0,b)$. Thus, $V_{0,b}^{E}(x)$ is concave on $(0,b)$.\\

\noindent(ii) Taking differentiation on \eqref{9823-1} yields
\begin{eqnarray}
V_{0,b^E}^{E\,\prime\prime}(x)
\hspace{-0.3cm}&=&\hspace{-0.3cm}
-\delta W_{\delta}(b^E-x)-\frac{\phi-Z_{\delta}(b^E)}{W_{\delta}(b^E)}W_{\delta}^\prime(b^E-x),\quad
 x\in(0,b^E).\label{9823-2}
\end{eqnarray}
By taking limits $x\uparrow b^E$, we have
\begin{eqnarray}
V_{0,b^E}^{E\,\prime\prime}(b^E-)
\hspace{-0.3cm}&=&\hspace{-0.3cm}
-\delta W_{\delta}(0)-\frac{\phi-Z_{\delta}(b^E)}{W_{\delta}(b^E)}W_{\delta}^\prime(0+)=0,\label{9823-3}
\end{eqnarray}
where the last equality follows by noting $\delta W_{\delta}(0)=0$ when $X$ is of unbounded variation (see \eqref{9823-4}) and $ \phi-Z_{\delta}(b^E)=0$ (see \eqref{9823-5}).
\hfill $\square$\\

\noindent
\textbf{Proof for Lemma \ref{Vb}}.
We derive explicit expression of $V_{0,b}(x)$ by finding representations for the two expected integral terms in \eqref{expdouble} separately utilising existing results in the literature. 
Note that
\begin{align}
 &\mathrm{E}_{x}\bigg[
\int^{\infty}_{0} e^{-(\delta +\lambda) t} \mathcal{P}^E(U^{\pi^{0,b}}{(t)};\pi^{0,b})dt
\bigg]\nonumber\\
=&\int^{b}_{0}\left( \int_0^{\infty} e^{-(\delta +\lambda) t} \mathcal{P}^E(y;\pi^{0,b})dt\right)\mathrm{P}_x(U^{\pi^{0,b}}{(t)}\in dy)\nonumber\\
=&\int^{b}_{0}\left( \int_0^{\infty} e^{-(\delta +\lambda) t} V^E_{0,b}(y)dt\right)\mathrm{P}_x(U^{\pi^{0,b}}{(t)}\in dy)\label{8823-2}\\
=&\int^{\infty}_{0} e^{-(\delta +\lambda) t} \left( \int_0^{b}V^E_{0,b}(y)\mathrm{P}_x(U^{\pi^{0,b}}{(t)}\in dy)\right) dt\nonumber\\
=&\int_{0}^{b}V^E_{0,b}(y)
\left(\frac{Z_{\delta+\lambda}(b-x) )W_{\delta+\lambda}^{\prime+}(y)}{(\delta+\lambda)W_{\delta+\lambda}(b)}-W_{\delta+\lambda}(y-x)\right) {d}y
\nonumber\\
&
+ V^E_{0,b}(0)
\left(\frac{Z_{\delta+\lambda}(b-x) W_{\delta+\lambda}(0)}{(\delta+\lambda)W_{\delta+\lambda}(b)}\right),\quad x\in[0,b],\label{8823-1}
\end{align}
where \eqref{8823-2} follows by noting $\mathcal{P}^E(y;\pi^{0,b})=V^E_{0,b}(y)$ (see \eqref{def-vobE}),  the  equality in \eqref{8823-1} follows by Theorem 1 of \cite{PM03} and $W_{\delta+\lambda}^{\prime+}$ there represents the right derivative of $W_{\delta+\lambda}$.
Noting that $\frac{Z_{\delta+\lambda}^\prime(\cdot)}{\delta+\lambda}=W_{\delta+\lambda}(\cdot)$ (see \eqref{8823-3}), it follows that
\begin{align}
&\int_{0}^{b}V^E_{0,b}(y)
\left(\frac{Z_{\delta+\lambda}(b-x) )W_{\delta+\lambda}^{\prime+}(y)}{(\delta+\lambda)W_{\delta+\lambda}(b)}-W_{\delta+\lambda}(y-x)\right) {d}y\nonumber\\
=& \int_{0}^{b}V^E_{0,b}(y)
\left(\frac{Z_{\delta+\lambda}(b-x) )W_{\delta+\lambda}^{\prime+}(y)}{(\delta+\lambda)W_{\delta+\lambda}(b)}-\frac{Z_{\delta+\lambda}^\prime(y-x)}{\delta+\lambda}\right) {d}y\nonumber\\
=&-V^{E}_{0,b}(0)\frac{Z_{\delta+\lambda}(b-x) )W_{\delta+\lambda}(0)}{(\delta+\lambda)W_{\delta+\lambda}(b)}
\nonumber\\
&
- \int_{0}^{b}V^{E\,\prime}_{0,b}(y)
\left(\frac{Z_{\delta+\lambda}(b-x) )W_{\delta+\lambda}(y)}{(\delta+\lambda)W_{\delta+\lambda}(b)}-\frac{Z_{\delta+\lambda}(y-x)}{\delta+\lambda}\right) {d}y,\quad x\in[0,b],\label{8823-4}
\end{align}
where the last equality is obtained by applying integration by parts. Now combining \eqref{8823-1} and \eqref{8823-4} yields that
\begin{align}
 &\mathrm{E}_{x}\bigg[
\int^{\infty}_{0} e^{-(\delta +\lambda) t} \mathcal{P}^E(U^{\pi^{0,b}}{(t)};\pi^{0,b})dt
\bigg]
\nonumber\\
=& \int_{0}^{b}V^{E\,\prime}_{0,b}(y)
\left(\frac{Z_{\delta+\lambda}(y-x)}{\delta+\lambda}-\frac{Z_{\delta+\lambda}(b-x) )W_{\delta+\lambda}(y)}{(\delta+\lambda)W_{\delta+\lambda}(b)}\right) {d}y, \quad x\in[0,b].
\label{8823-5}
\end{align}

In addition, mathematically, the sum of the first and second term on the right hand side of \eqref{def2-v0b} is the same as the expected payoff of the double-barrier strategy under exponential discounting with the discount rate $\delta+\lambda$, and hence the same as $\bar{v}_b$ in \cite{Bayraktar13} with $\delta$ there being replaced by $\delta+\lambda$ . Thus, by using results in  \cite{Bayraktar13} we have
\begin{align}
    &\mathrm{E}_{x}\bigg[\int^{\infty}_{0} e^{-(\delta+\lambda) t}(d L^b{(t)}-\phi dR^0{(t)})
\bigg]\nonumber\\
=&-
\overline{Z}_{\delta+\lambda}(b-x)
-\frac{\psi^{\prime}(0+)}{\delta+\lambda}
+\frac{Z_{\delta+\lambda}(b-x)}{\delta W_{\delta+\lambda}(b)}
\Big[Z_{\delta+\lambda}(b)-\phi
\Big], \quad x\in[0,b].\label{8823-6}
\end{align}
Furthermore, if we write $f(x):=\mathrm{E}_{x}\bigg[\int^{\infty}_{0} e^{-(\delta+\lambda) t}(d L^b{(t)}-\phi dR^0{(t)})
\bigg]$, then
\begin{align}\label{8823-7}
(\mathcal{A}-(\delta+\lambda))f(x)=0,\quad x\in(0,b).
\end{align}
Define $g(x):=\mathrm{E}_{x}\bigg[
\int^{\infty}_{0} e^{-(\delta +\lambda) t} \mathcal{P}^E(U^{\pi^{0,b}}{(t)};\pi^{0,b})dt
\bigg]= \int_{0}^{b}V^{E\,\prime}_{0,b}(y)
\left(\frac{Z_{\delta+\lambda}(y-x)}{\delta+\lambda}-\frac{Z_{\delta+\lambda}(b-x) )W_{\delta+\lambda}(y)}{(\delta+\lambda)W_{\delta+\lambda}(b)}\right) {d}y$. Then $V_{0,b}(x)=f(x)+g(x)$. 
Combining \eqref{def2-v0b}, \eqref{8823-5} and \eqref{8823-6} we can conclude that for $b\in(0,\infty)$,
\begin{eqnarray}
\label{v0b(x)}
\hspace{-0.5cm}
V_{0,b}(x)
\hspace{-0.3cm}&=&\hspace{-0.3cm}
\left\{
\small
\begin{aligned}
&
-
\overline{Z}_{\delta+\lambda}(b-x)
-\frac{\psi^{\prime}(0+)}{\delta+\lambda}
+
\frac{\lambda\beta}{{\delta+\lambda}}\bigg[V_{0,b}^E(0)+
\int_{0}^{b}V_{0,b}^{E\,\prime}(y)
Z_{\delta+\lambda}(y-x){d}y\bigg]&
\\&
+\frac{Z_{\delta+\lambda}(b-x)}{(\delta+\lambda)W_{\delta+\lambda}(b)}
\Big[Z_{\delta+\lambda}(b)-\phi
-\lambda\beta\int_{0}^{b}V_{0,b}^{E\,\prime}(y)W_{\delta+\lambda}(y)
{d}y\Big],& x\in[0,b],\\
&x-b+V_{0,b}(b),& x\in(b,\infty),\\
&\phi x+V_{0,b}(0), &x\in(-\infty,0),
\end{aligned}
\right.
\end{eqnarray}
where the expressions on $(b,\infty)$ and $(-\infty,0)$ follows directly from the construction of the double-barrier strategy.
Since the function $Z_{\delta+\lambda}(x)$ is continuous (hence, bounded) on $[-b,b]$, is continuously differentiable on $[-b,b]\setminus\{0\}$ and is left and right-differentiable at $x=0$. Hence, the function $Z_{\delta+\lambda}^{\prime}$ is continuous on $[-b,b]\setminus\{0\}$ and is bounded on $[-b,b]$, where, $Z_{\delta+\lambda}^{\prime}(0)$ is understood as the left or right-derivative at $0$.
Then, an application of the bounded convergence theorem yields
\begin{eqnarray}
\label{v0b'(x)}
V_{0,b}^{\prime}(x)
\hspace{-0.3cm}&=&\hspace{-0.3cm}
\left\{
\small
\begin{aligned}
&
Z_{\delta+\lambda}(b-x)-\lambda\beta\int_0^bV_{0,b}^{E\,\prime}(y)W_{\delta+\lambda}(y-x){d}y
&
\\&
-\frac{W_{\delta+\lambda}(b-x)}{W_{\delta+\lambda}(b)}\Big[Z_{\delta+\lambda}(b)-\phi-\lambda\beta\int_0^bV_{0,b}^{E\,\prime}(y)W_{\delta+\lambda}(y){d}y\Big],& x\in(0,b),\\
&1,& x\in(b,\infty),\\
&\phi, &x\in(-\infty,0),
\end{aligned}
\right.
\end{eqnarray}
which implies that $V_{0,b}(x)$ is continuously differentiable over $(-\infty,\infty)\setminus\{b\}$. And, furthermore, if $X$ has paths of unbounded variation, $V_{0,b}(x)$ can be checked to be twice continuously differentiable over $(0,\infty)\setminus\{b\}$ using the boundedness of $W^{\prime}_{\delta+\lambda}$ on $[-b,b]$ and bounded convergence theorem. Moreover,
if $X$ is of unbounded variation, $V_{0,b}^\prime(b)=1$ by letting $x\uparrow b$ on \eqref{v0b'(x)} and using  $Z_{\delta+\lambda}(0)=1$ (see \eqref{9823-9}) and $W_{\delta+\lambda}(0)=0$ (see \eqref{9823-4}).

Recall that the function $Z_{\delta+\lambda}(x)$ is continuous (hence, bounded) on $[-b,b]$, is continuously differentiable on $[-b,b]\setminus\{0\}$ and is left and right-differentiable at $x=0$. Hence, the function $Z_{\delta+\lambda}^{\prime}$ is continuous on $[-b,b]\setminus\{0\}$ and is bounded on $[-b,b]$, where, $Z_{\delta+\lambda}^{\prime}(0)$ is understood as the left or right-derivative at $0$. Moreover, if $\sigma>0$, then the function $Z_{\delta+\lambda}$ is twice continuously differentiable
on $[-b,b]$ (hence, $Z^{\prime\prime}_{\delta+\lambda}$ is bounded on $[-b,b]$). Therefore, 
using the bounded convergence theorem and the Fubini's theorem, one can obtain
\begin{eqnarray}
\hspace{-0.3cm}&&\hspace{-0.3cm}
(\mathcal{A}-({\delta+\lambda}))\int_0^bV_{0,b}^{E\,\prime}(y)Z_{\delta+\lambda}(y-x){d}y
=
\int_0^bV_{0,b}^{E\,\prime}(y)(\mathcal{A}-({\delta+\lambda}))Z_{\delta+\lambda}(y-x){d}y
\nonumber\\
\hspace{-0.3cm}&=&\hspace{-0.3cm}
\int_0^xV_{0,b}^{E\,\prime}(y)(\mathcal{A}-({\delta+\lambda}))Z_{\delta+\lambda}(y-x){d}y
+\int_x^bV_{0,b}^{E\,\prime}(y)(\mathcal{A}-({\delta+\lambda}))Z_{\delta+\lambda}(y-x){d}y
\nonumber\\
\hspace{-0.3cm}&=&\hspace{-0.3cm}-({\delta+\lambda})\int_0^xV_{0,b}^{E\,\prime}(y){d}y
\nonumber\\
\hspace{-0.3cm}&=&\hspace{-0.3cm}-({\delta+\lambda})\left(V_{0,b}^{E}(x)-V_{0,b}^{E}(0)\right),\quad x\in(0,b),\nonumber
\end{eqnarray}
where, in the third equality, we have also used the facts that $(\mathcal{A}-({\delta+\lambda}))Z_{\delta+\lambda}(y-x)=0$ for $x\in(0,y)$ (see \cite{Bayraktar13})
and $Z_{\delta+\lambda}(y-x)=1$ for $x\in(y,b)$.
Hence, by \eqref{v0b(x)}, we have
\begin{eqnarray}
\hspace{-0.3cm}&&\hspace{-0.3cm}
\mathcal{A}V_{0,b}(x)-({\delta+\lambda}) V_{0,b}(x)+\lambda\beta V_{0,b}^E(x)
\nonumber\\
\hspace{-0.3cm}&=&\hspace{-0.3cm}
(\mathcal{A}-{\delta+\lambda})\left(-\overline{Z}_{\delta+\lambda}(b-x)-\frac{\psi^{\prime}(0+)}{{\delta+\lambda}}+\frac{Z_{{\delta+\lambda}}(b-x)}{({\delta+\lambda})W_{{\delta+\lambda}}(b)}
\Big[Z_{{\delta+\lambda}}(b)-\phi
-\lambda\beta\int_{0}^{b}V_{0,b}^{E\,\prime}(y)W_{{\delta+\lambda}}(y)
{d}y\Big]\right)
\nonumber\\
\hspace{-0.3cm}&&\hspace{-0.3cm}
+(\mathcal{A}-({\delta+\lambda}))\frac{\lambda\beta}{{\delta+\lambda}}\left(V_{0,b}^E(0)+\int_0^bV_{0,b}^{E\,\prime}(y)Z_{\delta+\lambda}(y-x){d}y\right)+\lambda\beta V_{0,b}^E(x)
\nonumber\\
\hspace{-0.3cm}&=&\hspace{-0.3cm}0,\quad x\in(0,b).
\end{eqnarray}
The proof is complete.
\hfill $\square$\\

\noindent \textbf{Proof for Lemma \ref{lem.3.5}}. 
Note that $\ell(0)=1-\phi<0$ (due to $\phi>1$), and so  $b^{\ast}>0$. By \eqref{expdouble} and noticing $V^{E}_{0,b^E}(\cdot)=\mathcal{P}^E(\cdot;\pi^{0,b})$ (see \eqref{def-vobE}), it follows that
\begin{eqnarray}
V_{0,b^E}(x)
\hspace{-0.3cm}&=&\hspace{-0.3cm}
\mathrm{E}_x\bigg[\int_0^{\tau}e^{-\delta t}{d}L^{b^E}{(t)}-\phi\int_0^{\tau}e^{-\delta t}{d}R^{0}{(t)}+\beta e^{-\delta \tau} V^{E}_{0,b^E}(U^{0,b^E}{(\tau)})\bigg].\nonumber
\end{eqnarray}
By the strong Markov property we can obtain
\begin{eqnarray}
V_{0,b^E}^E(x)
\hspace{-0.3cm}&=&\hspace{-0.3cm}
\mathrm{E}_x\bigg[\int_0^{\infty}e^{-\delta t}dL^{b^E}{(t)}-\phi\int_0^{\infty}e^{-\delta t}{d}R^{0}{(t)}\bigg]
\nonumber\\
\hspace{-0.3cm}&=&\hspace{-0.3cm}
\mathrm{E}_x\bigg[\int_0^{\tau}e^{-\delta t}{d}L^{b^E}{(t)}-\phi\int_0^{\tau}e^{-\delta t}{d}R^{0}{(t)}+ e^{-\delta \tau} V^{E}_{0,b^E}(U^{0,b^E}{(\tau)})\bigg]
.
\end{eqnarray}
Taking difference yields
\begin{eqnarray}
V_{0,b^E}(b^E)-V_{0,b^E}^E(b^E)
\hspace{-0.3cm}&=&\hspace{-0.3cm}
-(1-\beta)\mathrm{E}_x\bigg[e^{-\delta \tau} V^{E}_{0,b^E}(U^{0,b^E}{(\tau)})\bigg]
\nonumber\\
\hspace{-0.3cm}&=&\hspace{-0.3cm}
-(1-\beta)\mathrm{E}_x\bigg[\int_0^{\infty}\lambda e^{-(\delta+\lambda) s} V^{E}_{0,b^E}(U^{0,b^E}{(s)}){d}s\bigg]
\nonumber\\
\hspace{-0.3cm}&\geq&\hspace{-0.3cm}
-(1-\beta)V_{0,b^E}^E(b^E)\mathrm{E}_x\bigg[\int_0^{\infty}\lambda e^{-(\delta+\lambda) s}{d}s\bigg]
\nonumber\\
\hspace{-0.3cm}&=&\hspace{-0.3cm}
-(1-\beta)\frac{\lambda}{\lambda+\delta}V_{0,b^E}^E(b^E),\nonumber
\end{eqnarray}
where the last inequality follows by the fact that $V_{0,b^E}^E(x)$ is  increasing on $[0,b^E)$ and $\beta\le 1$. The last equation implies
\begin{eqnarray}\label{Vb-VbE.2}
(\lambda+\delta)V_{0,b^E}(b^E)\geq(\lambda\beta+\delta)V_{0,b^E}^E(b^E).
\end{eqnarray}
Distinguish the following three mutually exclusive and collectively exhaustive cases.

\begin{itemize}
    \item[(a)]
When $\sigma\in(0,\infty)$ (hence, $X$ has paths of unbounded variation, and, $W_{q}^{\prime}(0+)=2/\sigma^{2}$), by \eqref{VE.hjb} (resp., \eqref{V.hjb}) and the twice continuous differentiability of $V_{0,b}^{E}(x)$ (resp., $V_{0,b}(x)$) on $(0,b)$ for $b\in(0,\infty)$, one has
\begin{eqnarray}\label{hjb.1.v1}
\hspace{-0.5cm}\frac{\sigma^2}{2}V_{0,b}^{E\,\prime\prime}(b-)
\hspace{-0.3cm}&=&\hspace{-0.3cm}\gamma V_{0,b}^{E\,\prime}(b-)-\int_0^{\infty}\left(V_{0,b}^{E}(b+y)-V_{0,b}^{E}(b)-V_{0,b}^{E\,\prime}(b-)y\mathbf{1}_{(0,1)}(y)\right)\nu({d}y)
\nonumber\\
\hspace{-0.3cm}&&\hspace{-0.3cm}
+\delta V_{0,b}^E(b),
\\
\label{hjb.2.v1}
\frac{\sigma^2}{2}V_{0,b}^{\prime\prime}(b-)
\hspace{-0.3cm}&=&\hspace{-0.3cm}\gamma V_{0,b}^{\prime}(b-)-\int_0^{\infty}\left(V_{0,b}(b+y)-V_{0,b}(b)-V_{0,b}^{\prime}(b-)y\mathbf{1}_{(0,1)}(y)\right)\nu({d}y)
\nonumber\\
\hspace{-0.3cm}&&\hspace{-0.3cm}
+(\delta+\lambda) V_{0,b}(b)-\lambda\beta V_{0,b}^E(b).
\end{eqnarray}
By setting  $b=b^E$  in \eqref{hjb.1.v1} and \eqref{hjb.2.v1}, taking difference of the two equation and then using  $V_{0,b^E}^{E\,\prime\prime}(b^E)=0$ (see Lemma \ref{lem2.4}(ii)) we can obtain
\begin{eqnarray}\label{Vb-VbE}
\frac{\sigma^2}{2}V_{0,b^E}^{\prime\prime}(b^E-)
\hspace{-0.3cm}&=&\hspace{-0.3cm}\gamma(V_{0,b^E}^{\prime}(b^E-)-V_{0,b^E}^{E\,\prime}(b^E-))+(\delta+\lambda)V_{0,b^E}(b^E)-(\lambda\beta+\delta)V_{0,b^E}^E(b^E)
\nonumber\\
\hspace{-0.3cm}&&\hspace{-0.3cm}
-\int_0^{\infty}\left(V_{0,b^E}(b^E+y)-V_{0,b^E}(b^E)-V_{0,b^E}^{\prime}(b^E-)y\mathbf{1}_{(0,1)}(y)\right)\nu({d}y)
\nonumber\\
\hspace{-0.3cm}&&\hspace{-0.3cm}
+\int_0^{\infty}\left(V_{0,b^E}^{E}(b^E+y)-V_{0,b^E}^{E}(b^E)-V_{0,b^E}^{E\,\prime}(b^E-)y\mathbf{1}_{(0,1)}(y)\right)\nu({d}y)
\nonumber\\
\hspace{-0.3cm}&=&\hspace{-0.3cm}
({\delta+\lambda})V_{0,b^E}(b^E)-(\lambda\beta+\delta)V_{0,b^E}^E(b^E),
\end{eqnarray}
where the last equation follows by noting $V_{0,b^E}^{E}(b^E+y)-V_{0,b^E}^{E}(b^E)=y=V_{0,b^E}(b^E+y)-V_{0,b^E}(b^E)$ (see \eqref{V.0bE.x} and \eqref{v0b(x)}) and  $V_{0,b^E}^{E\,\prime}(b^E-)=1=V_{0,b^E}^{\prime}(b^E-)$ (see \eqref{9823-7} and \eqref{9823-8}) when $X$ has paths of unbounded variation. On the other end, it follows by 
 \eqref{v0b(x)} that 
\begin{align}\label{condition.1.v1}
\frac{\sigma^2}{2}V_{0,b^E}^{\prime\prime}(b^E-)=\frac{W_{{\delta+\lambda}}^{\prime}(0+)}{W_{{\delta+\lambda}}(b^{E})}
\bigg[Z_{{\delta+\lambda}}(b^{E})-\phi-\lambda\beta\int_{0}^{b^{E}}V_{0,b^{E}}^{E\,\prime}(y)W_{{\delta+\lambda}}(y)
{d}y\bigg].
\end{align}
Therefore, combining \eqref{Vb-VbE} and \eqref{condition.1.v1} gives
\begin{align}
&\frac{W_{{\delta+\lambda}}^{\prime}(0+)}{W_{{\delta+\lambda}}(b^{E})}
\bigg[Z_{{\delta+\lambda}}(b^{E})-\phi-\lambda\beta\int_{0}^{b^{E}}V_{0,b^{E}}^{E\,\prime}(y)W_{{\delta+\lambda}}(y)
{d}y\bigg]\nonumber\\
=&
({\delta+\lambda})V_{0,b^E}(b^E)-(\lambda\beta+\delta)V_{0,b^E}^E(b^E)\geq0,\nonumber
\end{align}
where the last inequality follows by
 \eqref{Vb-VbE.2}. Hence, we can conclude that 
\begin{eqnarray}\label{condition.1}
Z_{{\delta+\lambda}}(b^{E})-\phi-\lambda\beta\int_{0}^{b^{E}}V_{0,b^{E}}^{E\,\prime}(y)W_{{\delta+\lambda}}(y)
{d}y
\hspace{-0.3cm}&\geq&\hspace{-0.3cm}0,
\end{eqnarray}
which along with the definition of $b{^\ast}$ implies $b^E\ge b{^\ast}$.
\item[(b)]
When $\sigma=0$ and $X$ has paths of unbounded variation (i.e., $\int_{0}^{1}y\nu({d}y)=\infty$), it follows by   \eqref{VE.hjb} and \eqref{V.hjb}) that
\begin{eqnarray}\label{hjb.1}
\hspace{-0.5cm}0
\hspace{-0.3cm}&=&\hspace{-0.3cm}\gamma V_{0,b^E}^{E\,\prime}(b^E-)-\int_0^{\infty}\left(V_{0,b^E}^{E}(b^E+y)-V_{0,b^E}^{E}(b^E)-V_{0,b^E}^{E\,\prime}(b^E-)y\mathbf{1}_{(0,1)}(y)\right)\nu({d}y)
\nonumber\\
\hspace{-0.3cm}&&\hspace{-0.3cm}
+\delta V_{0,b^E}^E(b^E),
\\
\label{hjb.2}
0
\hspace{-0.3cm}&=&\hspace{-0.3cm}\gamma V_{0,b^E}^{\prime}(b^E-)-\int_0^{\infty}\left(V_{0,b^E}(b^E+y)-V_{0,b^E}(b^E)-V_{0,b^E}^{\prime}(b^E-)y\mathbf{1}_{(0,1)}(y)\right)\nu({d}y)
\nonumber\\
\hspace{-0.3cm}&&\hspace{-0.3cm}
+(\delta+\lambda) V_{0,b^E}(b^E)-\lambda\beta V_{0,b^E}^E(b^E).
\end{eqnarray}
Taking differences leads to \begin{eqnarray}
0
\hspace{-0.3cm}&=&\hspace{-0.3cm}\gamma(V_{0,b^E}^{\prime}(b^E-)-V_{0,b^E}^{E\,\prime}(b^E-))+({\delta+\lambda})V_{0,b^E}(b^E)-(\lambda\beta+\delta)V_{0,b^E}^E(b^E)
\nonumber\\
\hspace{-0.3cm}&&\hspace{-0.3cm}
-\int_0^{\infty}\left(V_{0,b^E}(b^E+y)-V_{0,b^E}(b^E)-V_{0,b^E}^{\prime}(b^E-)y\mathbf{1}_{(0,1)}(y)\right)\nu({d}y)
\nonumber\\
\hspace{-0.3cm}&&\hspace{-0.3cm}
+\int_0^{\infty}\left(V_{0,b^E}^{E}(b^E+y)-V_{0,b^E}^{E}(b^E)-V_{0,b^E}^{E\,\prime}(b^E-)y\mathbf{1}_{(0,1)}(y)\right)\nu({d}y)
\label{2.20.eq-1}\\
\hspace{-0.3cm}&=&\hspace{-0.3cm}
(\delta+\lambda)V_{0,b^E}(b^E)-(\lambda\beta+\delta)V_{0,b^E}^E(b^E),\label{2.20.eq}
\end{eqnarray}
where the last equality follows by noting $V_{0,b^E}^{E}(b^E+y)-V_{0,b^E}^{E}(b^E)=y=V_{0,b^E}(b^E+y)-V_{0,b^E}(b^E)$ (by \eqref{V.0bE.x} and \eqref{v0b(x)}),  and  $V_{0,b^E}^{E\,\prime}(b^E-)=1=V_{0,b^E}^{\prime}(b^E-)$ (see \eqref{9823-7} 
and \eqref{9823-8}) as $X$ has paths of unbounded variation.
From  \eqref{V.0bE.x} and \eqref{v0b(x)} it follows
\begin{align}
V^E_{0,b^E}(b^E)
&=
-
\overline{Z}_{\delta}(0)
-\frac{\psi^{\prime}(0+)}{\delta}
+\frac{Z_{\delta}(0)}{\delta W_{\delta}(b^E)}
\Big[Z_{\delta}(b^E)-\phi
\Big]\nonumber\\
&=
-\frac{\psi^{\prime}(0+)}{\delta}
+\frac{1}{\delta W_{\delta}(b^E)}
\Big[Z_{\delta}(b^E)-\phi
\Big]\nonumber\\
&=-\frac{\psi^{\prime}(0+)}{\delta},\label{128-1}
\end{align}
where the last equality is due to $Z_{\delta}(b^E)-\phi=0$ (see \eqref{9823-5}),
and
\begin{align}
V_{0,b^E}(b^E)=&-
\overline{Z}_{\delta+\lambda}(0)
-\frac{\psi^{\prime}(0+)}{\delta+\lambda}
+
\frac{\lambda\beta}{{\delta+\lambda}}\bigg[V_{0,b^E}^E(0)+
\int_{0}^{b^E}V_{0,b^E}^{E\,\prime}(y)
Z_{\delta+\lambda}(y-b^E){d}y\bigg]
\nonumber\\
&+\frac{Z_{\delta+\lambda}(0)}{(\delta+\lambda)W_{\delta+\lambda}(b^E)}
\Big[Z_{\delta+\lambda}(b^E)-\phi
-\lambda\beta\int_{0}^{b^E}V_{0,b^E}^{E\,\prime}(y)W_{\delta+\lambda}(y)
{d}y\Big]\nonumber\\
=&
-\frac{\psi^{\prime}(0+)}{\delta+\lambda}
+
\frac{\lambda\beta}{{\delta+\lambda}}\bigg[V_{0,b^E}^E(0)+
\int_{0}^{b^E}V_{0,b^E}^{E\,\prime}(y)
{d}y\bigg]\nonumber\\
&
+\frac{1}{(\delta+\lambda)W_{\delta+\lambda}(b^E)}
\Big[Z_{\delta+\lambda}(b^E)-\phi
-\lambda\beta\int_{0}^{b^E}V_{0,b^E}^{E\,\prime}(y)W_{\delta+\lambda}(y)
{d}y\Big]\nonumber\\
=&
-\frac{\psi^{\prime}(0+)}{\delta+\lambda}
+
\frac{\lambda\beta}{{\delta+\lambda}}V_{0,b^E}^E(b^E)
+\frac{Z_{\delta+\lambda}(b^E)-\phi
-\lambda\beta\int_{0}^{b^E}V_{0,b^E}^{E\,\prime}(y)W_{\delta+\lambda}(y)
{d}y}{(\delta+\lambda)W_{\delta+\lambda}(b^E)},
\label{128-2}
\end{align}
where in the second to the last and the last equality we have used ${Z}_{\delta+\lambda}(y)=0$ and $\overline{Z}_{\delta+\lambda}(y)=0$ for $y\le 0$ (see \eqref{9823-9}), respectively. By plugging \eqref{128-1} and \eqref{128-2} into \eqref{2.20.eq} we arrive at
\begin{align}
\frac{1}{(\delta+\lambda)W_{\delta+\lambda}(b^E)}
\Big[Z_{\delta+\lambda}(b^E)-\phi
-\lambda\beta\int_{0}^{b^E}V_{0,b^E}^{E\,\prime}(y)W_{\delta+\lambda}(y)
{d}y\Big]=0,
\end{align}
which implies
\begin{eqnarray}
\label{condition.2}
Z_{{\delta+\lambda}}(b^{E})-\phi-\lambda\beta\int_{0}^{b^{E}}V_{0,b^{E}}^{E\,\prime}(y)W_{\delta+\lambda}(y){d}y=0.
\end{eqnarray}
This  along with the definition of $b{^\ast}$ implies $b^E\ge b{^\ast}$.

\item[(c)]
When $X$ has paths of bounded variation (i.e., $\sigma=0$ and 
$\int_{0}^{1}y\nu({d}y)<\infty$). Following the similar lines that lead to \eqref{2.20.eq-1}, we obtain the same result:
\begin{align}
0
=&\gamma(V_{0,b^E}^{\prime}(b^E-)-V_{0,b^E}^{E\,\prime}(b^E-))+(\delta+\lambda)V_{0,b^E}(b^E)-(\lambda\beta+\delta)V_{0,b^E}^E(b^E)
\nonumber\\
&
-\int_0^{\infty}\left(V_{0,b^E}(b^E+y)-V_{0,b^E}(b^E)-V_{0,b^E}^{\prime}(b^E-)y\mathbf{1}_{(0,1)}(y)\right)\nu({d}y)
\nonumber\\
=&
\left(\gamma+\int_0^1y\nu({d}y)\right)
(V_{0,b^E}^{\prime}(b^E-)-V_{0,b^E}^{E\,\prime}(b^E-))+(\delta+\lambda)V_{0,b^E}(b^E)-(\lambda\beta+\delta)V_{0,b^E}^E(b^E),\label{138-1}
\end{align}
where the last equality follows by noticing $V_{0,b^E}^{E}(b^E+y)-V_{0,b^E}^{E}(b^E)=y=V_{0,b^E}(b^E+y)-V_{0,b^E}(b^E)$ (by \eqref{V.0bE.x} and \eqref{v0b(x)}). 
Recall  from \eqref{9823-7} that  \begin{align}V_{0,b^E}^{E\,\prime}(b^E-)=1, \label{138-2}
    \end{align}  and note it follows from \eqref{v0b(x)} that 
\begin{align}
V_{0,b^E}^{\prime}(b^E-)=&Z_{\delta+\lambda}(0)-\lambda\beta\int_0^{b^E}V_{0,b^E}^{E\,\prime}(y)W_{\delta+\lambda}(y-b^E){d}y\nonumber\\
&
-\frac{W_{\delta+\lambda}(0)}{W_{\delta+\lambda}(b^E)}\Big[Z_{\delta+\lambda}(b^E)-\phi-\lambda\beta\int_0^{b^E}V_{0,b^E}^{E\,\prime}(y)W_{\delta+\lambda}(y){d}y\Big]\nonumber\\
=&1
-\frac{W_{\delta+\lambda}(0)}{W_{\delta+\lambda}(b^E)}\Big[Z_{\delta+\lambda}(b^E)-\phi-\lambda\beta\int_0^{b^E}V_{0,b^E}^{E\,\prime}(y)W_{\delta+\lambda}(y){d}y\Big],\label{138-3}
\end{align}
where the last equality follows by noting $Z_{\delta+\lambda}(0)=1$ and $Z_{\delta+\lambda}(z)=0$ for $z<0$.
Plugging \eqref{138-2} and \eqref{138-3}  into \eqref{138-1}
yields
\begin{eqnarray}
0\hspace{-0.3cm}&=&\hspace{-0.3cm}
-\left(\gamma+\int_0^1y\nu({d}y)\right)\frac{W_{\delta+\lambda}(0+)}{W_{\delta+\lambda}(b^E)}\bigg[Z_{{\delta+\lambda}}(b^{E})-\phi-\lambda\beta\int_{0}^{b^{E}}V_{0,b^{E}}^{E\,\prime}(y)W_{{\delta+\lambda}}(y)
{d}y\bigg]
\nonumber\\
\hspace{-0.3cm}&&\hspace{-0.3cm}
+({\delta+\lambda})V_{0,b^E}(b^E)-(\lambda\beta+\delta)V_{0,b^E}^E(b^E)
\nonumber\\
\hspace{-0.3cm}&\ge&\hspace{-0.3cm}
-\frac{Z_{{\delta+\lambda}}(b^{E})-\phi-\lambda\beta\int_{0}^{b^{E}}V_{0,b^{E}}^{E\,\prime}(y)W_{{\delta+\lambda}}(y)
{d}y}{W_{\delta+\lambda}(b^E)},\label{148-1}
\end{eqnarray}
where the second to the last equality follows by using $W_{\delta+\lambda}(0+)=\frac{1}{\gamma+\int_0^1y\nu({d}y)}$ (see \eqref{9823-4}) and  noting
$({\delta+\lambda})V_{0,b^E}(b^E)-(\lambda\beta+\delta)V_{0,b^E}^E(b^E)\ge 0$ (see 
 \eqref{Vb-VbE.2}).
 It follows immediately from the inequality \eqref{148-1} that \begin{eqnarray}
\label{condition.3}
Z_{\delta+\lambda}(b^{E})-\phi-\lambda\beta\int_{0}^{b^{E}}V_{0,b^{E}}^{E\,\prime}(y)W_{\delta+\lambda}(y){d}y\geq0,
\end{eqnarray}
which together the definition of $b^{\ast}$ (see \eqref{smoothcond.1}) implies $b^{\ast}\in(0, b^E]$. 
\end{itemize}
The proof is complete.
\hfill $\square$\\

\noindent
\textbf{Proof for Lemma \ref{b.ast}}
(a)\&(b). 
Let us first show the smoothness stated above. It follows by Lemma \ref{Vb} that 
$V_{0,b{^\ast}}(x)$
is continuously differentiable on $(-\infty,\infty)\setminus\{b\}$, and, if $X$ has paths of unbounded variation, $V_{0,b}(x)$ is continuously differentiable on $(-\infty,\infty)$ and twice continuously differentiable on $(0,\infty)\setminus\{b\}$. So we only need to prove the differentiabilty of $V_{0,b{^\ast}}$ at $0$ and $b{^\ast}$, and the twice differentiabilty at $b{^\ast}$ when $X$ has unbounded variation.
It follows by \eqref{v0b(x)-1}
that 
\begin{align}
V_{0,b{^\ast}}^{\prime}(x)
=&
Z_{\delta+\lambda}(b{^\ast}-x)-\lambda\beta\int_0^{b{^\ast}}V_{0,b{^\ast}}^{E\,\prime}(y)W_{\delta+\lambda}(y-x){d}y
\\
&
-\frac{W_{\delta+\lambda}(b{^\ast}-x)}{W_{\delta+\lambda}(b{^\ast})}\Big[Z_{\delta+\lambda}(b{^\ast})-\phi-\lambda\beta\int_0^{b{^\ast}}V_{0,b{^\ast}}^{E\,\prime}(y)W_{\delta+\lambda}(y){d}y\Big]\nonumber\\
=&Z_{\delta+\lambda}(b^{\ast}-x)
-\lambda\beta\int_{0}^{b^{\ast}}V_{0,b^{\ast}}^{E\,\prime}(y)W_{\delta+\lambda}(y-x){d}y
\label{14823-3}\\
=&Z_{\delta+\lambda}(b^{\ast}-x)
-\lambda\beta\int_{0}^{b^{\ast}}V_{0,b^{\ast}}^{E\,\prime}(y)W_{\delta+\lambda}(y-x){d}y
\nonumber\\
&+\bigg[\phi
+\lambda\beta\int_{0}^{b^{\ast}}V_{0,b^{\ast}}^{E\,\prime}(y)W_{\delta+\lambda}(y)
{d}y-Z_{\delta+\lambda}(b^{\ast})\bigg]\frac{Z_{\delta+\lambda}(b^{\ast}-x)}{Z_{\delta+\lambda}(b^{\ast})},\quad x\in(0,b{^\ast}),\label{14823-9}
\end{align}
 where the second to the last equality follows by noting $Z_{\delta+\lambda}(b{^\ast})-\phi-\lambda\beta\int_0^{b{^\ast}}V_{0,b{^\ast}}^{E\,\prime}(y)W_{\delta+\lambda}(y){d}y=0$ due to the  definition of $b^{\ast}$ (see \eqref{smoothcond.1}) and the fact that $0<b{^\ast}<\infty$ and the last equality follows by the same reason.  
From \eqref{14823-3} we can obtain
\begin{align}
&V_{0,b{^\ast}}^{\prime}(b{^\ast}-)
=
Z_{\delta+\lambda}(0)-\lambda\beta\int_0^{b{^\ast}}V_{0,b{^\ast}}^{E\,\prime}(y)W_{\delta+\lambda}(y-b{^\ast}){d}y=1=V_{0,b{^\ast}}^{\prime}(b{^\ast}+),
\label{14823-4}\\
&V_{0,b{^\ast}}^{\prime}(0+)
=Z_{\delta+\lambda}(b{^\ast})-\lambda\beta\int_0^{b{^\ast}}V_{0,b{^\ast}}^{E\,\prime}(y)W_{\delta+\lambda}(y){d}y=\phi=V_{0,b{^\ast}}^{\prime}(0-),
\label{14823-5}
\end{align}
where the second equality in \eqref{14823-4} and \eqref{14823-5} follow from  \eqref{v0b(x)-1} directly. The equations \eqref{14823-4} and \eqref{14823-5} imply that $V_{0,b{^\ast}}$ is also continuously differentiable at $b{^\ast}$ and $0$.
Now consider the situation that $X$ has paths of unbounded variation, by taking derivatives on  \eqref{14823-3} we can obtain
\begin{align}
&V_{0,b{^\ast}}^{\prime\prime}(b{^\ast}-)
=
W_{\delta+\lambda}(0)-\lambda\beta\int_0^{b{^\ast}}V_{0,b{^\ast}}^{E\,\prime}(y)W_{\delta+\lambda}^\prime(y-b{^\ast}){d}y=0=V_{0,b{^\ast}}^{\prime}(b{^\ast}+),
\end{align}
where the first equality in the above equation follows by noting $W_{\delta+\lambda}^\prime (z)=0$ for $z<0$ (due to \eqref{14823-7}) and $W_{\delta+\lambda}(0)=0$ (see \eqref{9823-4}), and the second equality follows directly from \eqref{Vb}. 
 
 We now proceed to show the concavity of $V_{0,b{^\ast}}$ on $(0,\infty)$. This can be achieved by proving that its derivative function is non-increasing.
Note that we already know that its derivative on $[b,\infty)$ and is a constant and thus non-increasing. So we only need to show the non-decreasing property of the function on $(0,b)$.  Define a surplus process $\{Y^{b^{\ast}}{(t)},t\geq 0\}$ as follows
\begin{eqnarray}
\label{def.V}
Y^{b^{\ast}}{(t)}
\hspace{-0.3cm}&:=&\hspace{-0.3cm}
X{(t)}-\sup_{0\leq s\leq t}(X{(s)}-b^{\ast})\vee 0,\quad t\geq 0,\nonumber
\end{eqnarray}
which is the spectrally positive L\'evy process $X$ with dividends deducted according to the barrier dividend strategy with barrier $b^{*}$.
We additionally define the bankruptcy time of $Y^{b^{\ast}}{(t)}$ as
\begin{eqnarray}
\label{def.zeta}
\zeta_{0}^{-}:=\inf\{t\geq0; Y^{b^{\ast}}{(t)}\leq 0\}.\nonumber
\end{eqnarray}
Then, by adapting Theorem 1 and Proposition 2 in \cite{Pistorius04}, one can get
\begin{align}
\label{Lap.zeta}
&\mathrm{E}_{x}\left[e^{-(\delta+\lambda)\zeta_{0}^{-}}\right]
=\frac{Z_{\delta+\lambda}(b^{\ast}-x)}{Z_{\delta+\lambda}(b^{\ast})},
\\
\label{Pot.Y}
&\int_{0}^{\infty}e^{-(\delta+\lambda)t}\mathrm{P}_{x}\left(Y^{b^{\ast}}{(t)}\in{d}y,\,t<\zeta_{0}^{-}\right){d}t
=
\bigg[\frac{Z_{\delta+\lambda}(b^{\ast}-x)}{Z_{\delta+\lambda}(b^{\ast})}W_{\delta+\lambda}(y)-W_{\delta+\lambda}(y-x)\bigg]\mathbf{1}_{[0,b^{\ast}]}(y){d}y.
\end{align}
It follows by \eqref{14823-9} that
\begin{align}
V_{0,b{^\ast}}^\prime(x)=&\phi \frac{Z_{\delta+\lambda}(b^{\ast}-x)}{Z_{\delta+\lambda}(b^{\ast})}
+\lambda\beta\int_{0}^{b^{\ast}}V_{0,b^{\ast}}^{E\,\prime}(y)\left[
\frac{Z_{\delta+\lambda}(b^{\ast}-x)}{Z_{\delta+\lambda}(b^{\ast})}W_{\delta+\lambda}(y)-W_{\delta+\lambda}(y-x)\right]
{d}y
\nonumber\\
=&\phi\mathrm{E}_{x}\left(e^{-(\delta+\lambda)\zeta_{0}^{-}}\right)+\lambda\beta
\int_{0}^{\infty}V_{0,b^{\ast}}^{E\,\prime}(y)\int_{0}^{\infty}e^{-(\delta+\lambda)t}\mathrm{P}_{x}\left(Y^{b^{\ast}}{(t)}\in{d}y,t<\zeta_{0}^{-}\right){d}t
\nonumber\\
=&\phi-
\mathrm{E}_{x}\left(\int_{0}^{\zeta_{0}^{-}}(\delta+\lambda)e^{-(\delta+\lambda)t}\left(\phi-\frac{\lambda\beta}{\delta+\lambda}V_{0,b^{\ast}}^{E\,\prime}(Y^{b^{\ast}}{(t)})\right){d}t\right),\quad x\in(0,b{^\ast}).\label{14823-10}
\end{align}
Note that both $Y^{b^{\ast}}{(t)}$ in non-decreasing in $X_0$ stochastically
and thus $\zeta_{0}^{-}$ is also non-decreasing in $X_0=x$ (see \cite{Pistorius04}). Recall $1\leq V_{0,b^{\ast}}^{E\,\prime}(x)\leq \phi$, $\lambda\ge 0$ and $\beta\le 1$. Hence, $ \phi-\frac{\lambda\beta}{\delta+\lambda}V_{0,b^{\ast}}^{E\,\prime}(Y^{b^{\ast}}{(t)})\ge 0$. Moreover, recall that 
$V_{0,b^{\ast}}^{E}(x)$ is concave (Lemma \ref{Vb.E} (ii)). Combining all the above together, we can see that   the function $$x\mapsto \phi-
\mathrm{E}_{x}\left[\int_{0}^{\zeta_{0}^{-}}(\delta+\lambda)e^{-(\delta+\lambda)t}\Big(\phi-\frac{\lambda\beta}{\delta+\lambda}V_{0,b^{\ast}}^{E\,\prime}(Y^{b^{\ast}}{(t)})\Big){d}t\right],\quad x\in[0,b^{\ast}],$$
is non-increasing, which together with \eqref{14823-10} implies that  $V_{0,b^{\ast}}^\prime(x)$ is non-decreasing on $(0,b^{\ast})$. This completes the proof.

\noindent (c) It follows immediately form \eqref{14823-4} that $V^{\prime}_{{0,b^\ast}}(b^\ast)=1$. 
\hfill $\square$\\

\noindent \textbf{Proof for Theorem \ref{16923-1}}.
Note by \eqref{expdouble} that
$V_{0,b^{\ast}}(x)=\mathcal{P}(x;\pi^{0,b^{\ast}},\pi^{0,b^{\ast}})\leq \sup_{\pi\in\Pi}\mathcal{P}(x;\pi,\pi^{0,b^{\ast}})$. It is sufficient to show that $V_{0,b^{\ast}}(x)\geq \sup_{\pi\in\Pi}\mathcal{P}(x;\pi,\pi^{0,b^{\ast}})$. Recall from Lemma \ref{b.ast} we know that $V_{0,b^\ast}(x)$ is continuously differentiable on $[0,\infty)$ and twice continuously differentiable on $[0,b^\ast)\cup(b^\ast,\infty)$, and if  $X$ has paths of unbounded variation, $V_{0,b^{\ast}}(x)$ is twice continuously differentiable on $[0,\infty)$. 

Consider an arbitrary admissible strategy, $\pi=(R,L)$. 
Note that $V_{0,b^\ast}$ is twice continuously differentiable on $(0,\infty)$.
By Theorem 4.57 (It\^{o}'s formula) in \cite{JaSh2003}, we can obtain that for $x\in(0,\infty)$,
\begin{eqnarray}
\hspace{-0.3cm}&&\hspace{-0.3cm}
e^{-(\delta+\lambda)t}V_{0,b^{\ast}}(U^{R,L}{(t)})-V_{0,b^{\ast}}(U^{R,L}({0-}))
\nonumber\\
\hspace{-0.3cm}&=&\hspace{-0.3cm}
-\int_{0-}^{t}(\delta+\lambda) e^{-(\delta+\lambda) s}V_{0,b^{\ast}}( {U}^{R,L}{(s-)}){d}s+\int_{0-}^{t} e^{-(\delta+\lambda) s}V_{0,b^{\ast}}^{\prime}( {U}^{R,L}{(s-)}){d} {U}^{R,L}{(s)}\nonumber\\
    \hspace{-0.3cm}&&\hspace{-0.3cm}
+\frac{1}{2} \int_{0-}^{t} e^{-(\delta+\lambda) s}V_{0,b^{\ast}}^{\prime\prime}( {U}^{R,L}{(s-)}){d}\langle  {U}^{R,L}(\cdot), {U}^{R,L}(\cdot)\rangle{(s)}\nonumber\\
    \hspace{-0.3cm}&&\hspace{-0.3cm}
+\sum_{s\leq t}e^{-(\delta+\lambda) s}\big(V_{0,b^{\ast}}( {U}^{R,L}{(s-)}+\Delta  {U}^{R,L}{(s)})-V_{0,b^{\ast}}( {U}^{R,L}{(s-)}))\Delta  {U}^{R,L}{(s)}\big)
\nonumber\\
\hspace{-0.3cm}&=&\hspace{-0.3cm}
-\int_{0-}^{t}(\delta+\lambda) e^{-(\delta+\lambda) s}V_{0,b^{\ast}}( {U}^{R,L}{(s-)}){d}s+\frac{1}{2} \sigma^2\int_{0-}^{t} e^{-(\delta+\lambda) s}V_{0,b^{\ast}}^{\prime\prime}( {U}^{R,L}{(s-)}){d}s
\nonumber\\
    \hspace{-0.3cm}&&\hspace{-0.3cm}
+\int_{0-}^{t} e^{-(\delta+\lambda) s}V_{0,b^{\ast}}^{\prime}( {U}^{R,L}{(s-)}){d} (-\gamma s +\sigma B(s)+M_1(s)+R^c{(s)}-L^c{(s)})\nonumber\\
 \hspace{-0.3cm}&&\hspace{-0.3cm}+\sum_{s\leq t}e^{-(\delta+\lambda) s}\big[V_{0,b^{\ast}}( {U}^{R,L}{(s-)}+\Delta X{(s)})-V_{0,b^{\ast}}( {U}^{R,L}{(s-)})
\big]\mathbf{1}_{\{\Delta X{(s)}\ge 1\}}
\nonumber\\
 \hspace{-0.3cm}&&\hspace{-0.3cm}
+\sum_{s\leq t}e^{-(\delta+\lambda) s}\big[V_{0,b^{\ast}}( {U}^{R,L}{(s-)}+\Delta M_1{(s)})-V_{0,b^{\ast}}( {U}^{R,L}{(s-)})
-V_{0,b^{\ast}}^{\prime}( {U}^{R,L}{(s-)})\Delta M_1{(s)}
\big]\mathbf{1}_{\{\Delta M_1{(s)}>0\}}
\nonumber\\
\hspace{-0.3cm}&&\hspace{-0.3cm}
+\sum_{s\leq t}e^{-(\delta+\lambda) s}\big[V_{0,b^{\ast}}( {U}^{R,L}{(s-)}+\Delta  {}R{(s)})-V_{0,b^{\ast}}( {U}^{R,L}{(s-)})
\big]\mathbf{1}_{\{\Delta R{(s)}>0\}}
\nonumber\\
\hspace{-0.3cm}&&\hspace{-0.3cm}
+\sum_{s\leq t}e^{-(\delta+\lambda) s}\big[V_{0,b^{\ast}}( {U}^{R,L}{(s-)}-\Delta  {}L{(s)})-V_{0,b^{\ast}}( {U}^{R,L}{(s-)})
\big]\mathbf{1}_{\{\Delta L{(s)}>0\}},\label{23823-2}
\end{eqnarray}
where $R^c$ and $D^c$ represent the continuous parts of $R$ and $D$ respectively, $\Delta{U}^{R,L}{(s)} = {U}^{R,L}{(s)} - {U}^{R,L}{(s-)}$, $\Delta{L}{(s)} = {L}{(s)} - {L}{(s-)}$, $\Delta X{(s)} = X{(s)} - X{(s-)}$, and $\Delta{M_1}{(s)} = {M_1}{(s)} - {M_1}{(s-)}$. The last equality follows by using the decomposition \eqref{decompostion}, and 
 noticing $X_1(t)$ has no continuous part. Notice that
 \begin{align}
 &\mathrm{E}_x\left[-\int_{0-}^{t}(\delta+\lambda) e^{-(\delta+\lambda) s}V_{0,b^{\ast}}( {U}^{R,L}{(s-)}){d}s+\frac{1}{2} \sigma^2\int_{0-}^{t} e^{-(\delta+\lambda) s}V_{0,b^{\ast}}^{\prime\prime}( {U}^{R,L}{(s-)}){d}s\right.\nonumber\\
 &+\left.\int_{0-}^{t} e^{-(\delta+\lambda) s}V_{0,b^{\ast}}^{\prime}( {U}^{R,L}{(s-)}){d} (-\gamma s +\sigma B(s)+M_1(s)+R^c{(s)}-L^c{(s)})\right]\nonumber\\ =&\mathrm{E}_x\left[\int_{0-}^{t} e^{-(\delta+\lambda) s}(\mathcal{A}-(\delta+\lambda))V_{0,b^{\ast}}( {U}^{R,L}{(s-)}){d}s\right.\nonumber\\
 &\left.-\int_{0-}^{t} e^{-(\delta+\lambda) s}\left(\int_0^\infty \left(V_{0,b^{\ast}}({U}{(s-)}+y)-V_{0,b^{\ast}}({U}^{R,L}{(s-)})-{V^{\prime}_{0,b^{\ast}}}({U}^{R,L}{(s-)})y\mathbf{1}\{|y|<1\}\right)v(dy)\right)
{d}s\right]\nonumber\\
&+\mathrm{E}_x\left[\int_{0-}^{t} e^{-(\delta+\lambda) s}V_{0,b^{\ast}}^{\prime}( {U}^{R,L}{(s-)}){d} (\sigma B(s)+M_1(s))\right]+\mathrm{E}_x\left[\int_{0-}^{t} e^{-(\delta+\lambda) s}V_{0,b^{\ast}}^{\prime}( {U}^{R,L}{(s-)}){d}R^c{(s)}\right]\nonumber\\
&-\mathrm{E}_x\left[\int_{0-}^{t} e^{-(\delta+\lambda) s}V_{0,b^{\ast}}^{\prime}( {U}^{R,L}{(s-)}){d}L^c{(s)}\right]
\nonumber\\
=&\mathrm{E}_x\left[\int_{0-}^{t} e^{-(\delta+\lambda) s}g_{b^{\ast}}( {U}^{R,L}{(s-)}){ds}\right]-\mathrm{E}_x\left[\int_{0-}^{t} e^{-(\delta+\lambda) s}\lambda \beta V_{0,b^{\ast}}^E( {U}^{R,L}{(s-)}){d}s\right]\nonumber\\
&-\mathrm{E}_x\left[\int_{0-}^{t} e^{-(\delta+\lambda) s}\left(\int_0^\infty \left(V_{0,b^{\ast}}({U}^{R,L}{(s-)}+y)-V_{0,b^{\ast}}({U}^{R,L}{(s-)})\right.\right.\right.\nonumber\\
&\left.\left.\left.-{V^{\prime}_{0,b^{\ast}}}({U}^{R,L}{(s-)})y\mathbf{1}\{|y|<1\}\right)v({d}y)\right)
{d}s\right]\nonumber\\
&+\mathrm{E}_x\left[M_2(t)\right]+\mathrm{E}_x\left[\int_{0-}^{t} e^{-(\delta+\lambda) s}V_{0,b^{\ast}}^{\prime}( {U}^{R,L}{(s-)}){d}R^c{(s)}-\int_{0-}^{t} e^{-(\delta+\lambda) s}V_{0,b^{\ast}}^{\prime}( {U}^{R,L}{(s-)}){d}L^c{(s)}\right], \label{23823-1}
\end{align}
where
\begin{align}\label{gb(x).def}
g_{b^\ast}(x)
:=&
(\mathcal{A}-(\delta+\lambda))V_{0,b^\ast}(x)+\lambda\beta V_{0,b^\ast}^E(x)
\nonumber \\
=&  
\frac{\sigma^2}{2}V_{0,b^\ast}^{\prime\prime}(x)-\gamma V_{0,b^\ast}^{\prime}(x)+\int_0^{\infty}\left(V_{0,b^\ast}(x+y)-V_{0,b^\ast}(x)-V_{0,b^\ast}^{\prime}(x)y\mathbf{1}_{(0,1)}(y)\right)\nu({d}y)
\nonumber\\
&-qV_{0,b^\ast}(x)+\lambda\beta V_{0,b^\ast}^E(x),\\
M_2(t):=&\int_{0-}^{t} e^{-(\delta+\lambda) s}V_{0,b^{\ast}}^{\prime}( {U}^{R,L}{(s-)}){d} (\sigma B(s)+M_1(s)).
\end{align}
As $\sigma B(t)+M_1(t)$ is a martingale and so $M_2(t)$ is a local martingale.
Further define 
\begin{align}
M_3(t):=&\sum_{s\leq t}e^{-(\delta+\lambda) s}\big[V_{0,b^{\ast}}( {U}^{R,L}{(s-)}+\Delta X{(s)})-V_{0,b^{\ast}}( {U}^{R,L}{(s-)})
\big]\mathbf{1}_{\{\Delta X{(s)}\ge 1\}}
\nonumber\\
 &
+\sum_{s\leq t}e^{-(\delta+\lambda) s}\times\nonumber\\
&\big[V_{0,b^{\ast}}( {U}^{R,L}{(s-)}+\Delta M_1{(s)})-V_{0,b^{\ast}}( {U}{(s-)})
-V_{0,b^{\ast}}^{\prime}( {U}^{R,L}{(s-)})\Delta M_1{(s)}
\big]\mathbf{1}_{\{\Delta M_1{(s)}>0\}}\nonumber\\
&-\int_{0-}^{t} e^{-(\delta+\lambda) s}\times\nonumber\\
&\left(\int_0^\infty \left(V_{0,b^{\ast}}({U}^{R,L}{(s-)}+y)-V_{0,b^{\ast}}({U}^{R,L}{(s-)})-{V^\prime_{0,b^{\ast}}}({U}^{R,L}{(s-)})y\mathbf{1}_{\{|y|<1\}}\right)v({d}y)\right)
{d}s.\label{23823-3}
\end{align}
Recall that both $M_1$ and $X$ have positive jumps only and the jumps of $M_1$ are identical to the jumps of the L\'evy process $X$ with magnitude smaller than unity and thus $\Delta M_1(s)=\Delta X(s)\mathbf{1}_{\{\Delta |X(s)|<1\}}$ and $\mathbf{1}_{\{\Delta M_1(s)>0\}}=\mathbf{1}_{\{0<\Delta X(s)<1\}}$. We can rewrite  $M_3(t)$ as
\begin{align*}
M_3(t):=&\sum_{s\leq t}e^{-(\delta+\lambda) s}\Big[\big(V_{0,b^{\ast}}( {U}^{R,L}{(s-)}+\Delta X{(s)})-V_{0,b^{\ast}}( {U}^{R,L}{(s-)})
\big)\mathbf{1}_{\{\Delta X{(s)}>0\}}
\nonumber\\
&\hspace{2cm}-V_{0,b^{\ast}}^{\prime}( {U}^{R,L}{(s-)})\Delta M_1{(s)}\mathbf{1}_{\{\Delta M_1{(s)}>0\}}\Big]
\nonumber\\
 &-\int_{0-}^{t} e^{-(\delta+\lambda) s}\left(\int_0^\infty \left(V_{0,b^{\ast}}({U}^{R,L}{(s-)}+y)-V_{0,b^{\ast}}({U}^{R,L}{(s-)})\right.\right.
 \nonumber\\
&\hspace{2.5cm}
\left.\left.-{V^\prime_{0,b^{\ast}}}({U}^{R,L}{(s-)})y\mathbf{1}_{\{|y|<1\}}\right)v({d}y)\right)
{d}s. 
\end{align*}
From the above we can see that $M_3(t)$ is a local martingale. 
Combining \eqref{23823-2}, \eqref{23823-1}
and \eqref{23823-3} yields
\begin{align}
&\mathrm{E}_x\left[e^{-qt}V_{0,b^{\ast}}(U^{R,L}{(t)})-V_{0,b^{\ast}}(U^{R,L}({0-}))\right]
\nonumber\\
=&\mathrm{E}_x\left[\int_{0-}^{t} e^{-(\delta+\lambda) s}g_{b^{\ast}}( {U}^{R,L}{(s-)}){ds}\right]-\mathrm{E}_x\left[\int_{0-}^{t} e^{-(\delta+\lambda) s}\lambda \beta V_{0,b^{\ast}}^E( {U}^{R,L}{(s-)})){d}s\right]\nonumber\\
&+\mathrm{E}_x\left[\int_{0-}^{t} e^{-(\delta+\lambda) s}V_{0,b^{\ast}}^{\prime}( {U}^{R,L}{(s-)}){d}R^c{(s)}-\int_{0-}^{t} e^{-(\delta+\lambda) s}V_{0,b^{\ast}}^{\prime}( {U}^{R,L}{(s-)}){d}L^c{(s)})\right]\nonumber\\
&+\mathrm{E}_x\left[M_2(t)\right]+\mathrm{E}_x\left[M_3(t)\right]\nonumber\\
&+\mathrm{E}_x\left[\sum_{s\leq t}e^{-(\delta+\lambda) s}\big[V_{0,b^{\ast}}( {U}^{R,L}{(s-)}+\Delta  {}R{(s)})-V_{0,b^{\ast}}( {U}^{R,L}{(s-)})
\big]\mathbf{1}_{\{\Delta R{(s)}>0\}}\right]
\nonumber\\
&
+\mathrm{E}_x\left[\sum_{s\leq t}e^{-(\delta+\lambda) s}\big[V_{0,b^{\ast}}( {U}^{R,L}{(s-)}-\Delta  {}L{(s)})-V_{0,b^{\ast}}( {U}^{R,L}{(s-)})
\big]\mathbf{1}_{\{\Delta L{(s)}>0\}}\right].
\label{23823-4}
\end{align}
 Recall that $M_2(t)$ and $M_3(t)$
are zero-mean local martingales. We can find a sequence of stopping times $(T_m)_{m\ge 1}$ with $\lim_{m\rightarrow \infty}T_m=\infty$ such that for each $m$, both $\{M_2(t\wedge T_m); {t\ge 0}\}$ and $\{M_3(t\wedge T_m);{t\ge 0}\}$ are martingales. As a result, for each $m$,
\begin{align}
\mathrm{E}_x\left[M_2(t\wedge T_m)\right]=0, \quad \mathrm{E}_x\left[M_3(t\wedge T_m)\right]=0.\label{23823-7}
\end{align}

It follows immediately by Lemma \ref{Vb} that
\begin{eqnarray}\label{g.b(x)=0}
g_{b^\ast}(x)=0,\quad x\in(0,b^\ast).
\end{eqnarray}
We now show that 
\begin{eqnarray}\label{g.b(x)<0}
g_{b^{\ast}}(x)\leq0,\quad x\in(0,\infty).
\end{eqnarray}
To this end, it is sufficient to show that $\limsup_{y\downarrow x}\frac{g_{b^{\ast}}(y)-g_{b^{\ast}}(x)}{y-x}\le 0$ for $x\ge b{^\ast}$. 
For $x>0$, let $y_{x,n}^-$ and $y_{x,n}^+$ represent the sequences with $y_{x,n}^-\uparrow x$ and $y_{x,n}^+\downarrow x$ as $
n\rightarrow\infty$ such that $$\limsup_{y\uparrow x}\frac{g_{b^{\ast}}(y)-g_{b^{\ast}}(x)}{y-x}=\lim_{n\rightarrow \infty}\frac{g_{b^{\ast}}(y_{x,n}^-)-g_{b^{\ast}}(x)}{y_{x,n}^--x}, \quad \limsup_{y\downarrow x}\frac{g_{b^{\ast}}(y)-g_{b^{\ast}}(x)}{y-x}=\lim_{n\rightarrow \infty}\frac{g_{b^{\ast}}(y_{x,n}^+)-g_{b^{\ast}}(x)}{y_{x,n}^+-x}.$$ 
It is now sufficient to show 
\begin{align}\lim_{n\rightarrow \infty}\frac{g_{b^{\ast}}(y_{x,n}^+)-g_{b^{\ast}}(x)}{y_{x,n}^+-x}\ge 0, \quad x\ge b{^\ast}. \label{22823-4}
\end{align}
First, let us consider the case of unbounded variation. Recall that in this case, $V_{0,b^{\ast}}$ and $V_{0,b^{\ast}}^E$ are both twice continuously differentiable on $(0,b^\ast)\vee (b^\ast, \infty)$.  We can find sub-sequences of $y_{x,n}^-$ and $y_{x,n}^+$, respectively, say $y_{x,n_k}^-$ and $y_{x,n_k}^+$, such that $\lim_{k\rightarrow\infty}\frac{V_{0,b^{\ast}}^{\prime\prime}(y_{x,n_k}^-)-V_{0,b^{\ast}}^{\prime\prime}(x)}{y_{x,n_k}^--x}$ and $\lim_{k\rightarrow\infty}\frac{V_{0,b^{\ast}}^{\prime\prime}(y_{x,n_k}^+)-V_{0,b^{\ast}}^{\prime\prime}(x)}{y_{x,n_k}^+-x}$  both exist. It follows by \eqref{g.b(x)=0} that $\lim_{n\rightarrow \infty}\frac{g_{b^{\ast}}(y_{x,n}^-)-g_{b^{\ast}}(x)}{y_{x,n}^--x}=0$.
Thus,  
\begin{eqnarray}\label{g(y-)-g(x)}
0
\hspace{-0.3cm}&=&\hspace{-0.3cm}
\lim_{k\rightarrow\infty}\frac{g_{b^{\ast}}(y_{x,n_k}^-)-g_{b^{\ast}}(x)}{y_{x,n_k}^--x}
\nonumber\\
\hspace{-0.3cm}&=&\hspace{-0.3cm}
\frac{\sigma^2}{2}\lim_{k\rightarrow\infty}\frac{V_{0,b^{\ast}}^{\prime\prime}(y_{x,n_k}^-)-V_{0,b^{\ast}}^{\prime\prime}(x)}{y_{x,n_k}^--x}-\gamma V_{0,b^{\ast}}^{\prime\prime}(x)
\nonumber\\
\hspace{-0.3cm}&&\hspace{-0.3cm}
+\int_0^{\infty}\left(V_{0,b^{\ast}}^{\prime}(x+y)-V_{0,b^{\ast}}^{\prime}(x)-V_{0,b^{\ast}}^{\prime\prime}(x)y\mathbf{1}_{(0,1)}(y)\right)\nu({d}y)
\nonumber\\
\hspace{-0.3cm}&&\hspace{-0.3cm}
-qV_{0,b^{\ast}}^{\prime}(x)+\lambda\beta V_{0,b^{\ast}}^{E\,\prime}(x),\quad x\in(0,b^{\ast}],
\end{eqnarray}
where the last equality follows by \eqref{gb(x).def}. 
Recall that $V_{0,b^{\ast}}(x)$ is concave and twice continuously differentiable (Lemma \ref{b.ast}). Hence, $V_{0,b^{\ast}}^{\prime\prime}(x)\leq0$ for $x\in(0,b^{\ast})$, which along with $V_{0,b^{\ast}}^{\prime\prime}(x)=0$ for $x\geq b^{\ast}$ (by \eqref{Vb}), yields
\begin{eqnarray}\label{V(3)}
\lim_{k\rightarrow\infty}\frac{V_{0,b^{\ast}}^{\prime\prime}(y^-_{b^{\ast},n_k})-V_{0,b^{\ast}}^{\prime\prime}(b^{\ast})}{y^-_{b^{\ast},n_k}-b^{\ast}}\geq0.
\end{eqnarray}
It follows by \eqref{gb(x).def} and by noting $V_{0,b^{\ast}}^{\prime}(x)=1$ for $x\geq b^{\ast}$ (by \eqref{Vb}) that 
$$g_{b^{\ast}}(x)=-\gamma +\int_0^{\infty}(y-y\mathbf{1}_{(0,1)}(y))\nu({d}y)-qV_{0,b^{\ast}}(x)+\lambda\beta V_{0,b^{\ast}}^E(x),\quad x\in[b^{\ast},\infty).$$  Thus,
\begin{eqnarray}
\lim_{k\rightarrow\infty}\frac{g_{b^{\ast}}(y_{x,n_k}^+)-g_{b^{\ast}}(x)}{y_{x,n_k}^+-x}
\hspace{-0.3cm}&=&\hspace{-0.3cm}-q V_{0,b^{\ast}}^{\prime}(x)+\lambda\beta V_{0,b^{\ast}}^{E\,\prime}(x)
\leq
-q+\lambda\beta V_{0,b^{\ast}}^{E\,\prime}(b^{\ast}),\quad x\in[b^{\ast},\infty),\label{22823-1}
\end{eqnarray}
 where the last inequality follows by using $V_{0,b^{\ast}}^\prime(x)=1$ for $x\geq b^{\ast}$ and the concavity of $V_{0,b^{\ast}}(x)$ again. 
Furthermore, setting $x=b^{\ast}$ in \eqref{g(y-)-g(x)}, and then using $V_{0,b^{\ast}}^{\prime}(x)=1$ for $x\geq b^{\ast}$ and $V_{0,b^{\ast}}^{\prime\prime}(x)=0$ for $x\geq b^{\ast}$ again we obtain
\begin{eqnarray}
0
\hspace{-0.3cm}&=&\hspace{-0.3cm}
\lim_{k\rightarrow\infty}\frac{g_{b^{\ast}}(y_{b^{\ast},n_k}^-)-g_{b^{\ast}}(b^{\ast})}{y_{b^{\ast},n_k}^--x}
\nonumber\\
\hspace{-0.3cm}&=&\hspace{-0.3cm}
\frac{\sigma^2}{2}\lim_{k\rightarrow\infty}\frac{V_{0,b^{\ast}}^{\prime\prime}(y_{b^{\ast},n_k}^-)-V_{0,b^{\ast}}^{\prime\prime}(b^{\ast})}{y_{b^{\ast},n_k}^--x}
-q+\lambda\beta V_{0,b^{\ast}}^{E\,\prime}(b^{\ast})
\nonumber\\
\hspace{-0.3cm}&\geq&\hspace{-0.3cm}
-q+\lambda\beta V_{0,b^{\ast}}^{E\,\prime}(x)
\nonumber\\
\hspace{-0.3cm}&=&\hspace{-0.3cm}
\lim_{k\rightarrow\infty}\frac{g_{b^{\ast}}(y_{x,n_k}^+)-g_{b^{\ast}}(x)}{y_{x,n_k}^+-x},\quad x\in[b^{\ast},\infty),\label{22823-5}
\end{eqnarray}
where the last inequality follows by \eqref{V(3)} and the last equality by \eqref{22823-1}.
Since $y_{x,n_k}^+$ is a sub-sequence of $y_{x,n}^+$,  \eqref{22823-5} implies that \eqref{22823-4} holds. 

Now, let us consider the case where $X$ has paths of bounded variation. In this case,  $\sigma=0$ in \eqref{gb(x).def}, and then
\begin{eqnarray}\label{gb(y)-gb(b).1}
0
\hspace{-0.3cm}&=&\hspace{-0.3cm}
\lim_{k\rightarrow\infty}\frac{g_{b^{\ast}}(y_{b^{\ast},n_k}^-)-g_{b^{\ast}}(b^{\ast})}{y_{b^{\ast},n_k}^--b^{\ast}}
\nonumber\\
\hspace{-0.3cm}&=&\hspace{-0.3cm}
-\gamma \frac{V_{0,b^{\ast}}^{\prime}(y_{b^{\ast},n_k}^-)-V_{0,b^{\ast}}^{\prime}(b^{\ast})}{y_{b^{\ast},n_k}^--b^{\ast}}-qV_{0,b^{\ast}}^{\prime}(b^{\ast})+\lambda\beta V_{0,b^{\ast}}^{E\,\prime}(b^{\ast})
\nonumber\\
\hspace{-0.3cm}&&\hspace{-0.3cm}
+\int_0^{\infty}\left(V_{0,b^{\ast}}^{\prime}(b^{\ast}+y)-V_{0,b^{\ast}}^{\prime}(b^{\ast})-\frac{V_{0,b^{\ast}}^{\prime}(y_{b^{\ast},n_k}^-)-V_{0,b^{\ast}}^{\prime}(b^{\ast})}{y_{b^{\ast},n_k}^--b^{\ast}}y\mathbf{1}_{(0,1)}(y)\right)\nu({d}y)
\nonumber\\
\hspace{-0.3cm}&=&\hspace{-0.3cm}
-q+\lambda\beta V_{0,b^{\ast}}^{E\,\prime}(b^{\ast})=-q+\lambda\beta ,
\end{eqnarray}
where  the last equality follows by noting that the left derivative of  $V_{0,b^{\ast}}^\prime(x)$ at $x=b^\ast$ is $1$ (Lemma \ref{b.ast}(c)), and $V_{0,b^{\ast}}^{E\,\prime}(b^{\ast})=1$ (by \eqref{V.0bE.x}). Furthermore, from \eqref{gb(x).def} and by noting $\sigma=0$ we know that for $x\ge b^\ast$, 
\begin{eqnarray}\label{gb(y)-gb(b).2}
\lim_{k\rightarrow\infty}\frac{g_{b^{\ast}}(y_{x,n_k}^+)-g_{b^{\ast}}(x)}{y_{x,n_k}^+-x}
\hspace{-0.3cm}&=&\hspace{-0.3cm}-q+\lambda\beta V_{0,b^{\ast}}^{E\,\prime+}(x)
=
-q+\lambda\beta,\quad x\in[b^{\ast},\infty),
\end{eqnarray}
where $V_{0,b^{\ast}}^{E\,\prime+}$ represents the right derivative of $V_{0,b^{\ast}}^E$ and  the last inequality follows by noting $V_{0,b^{\ast}}^{E\,\prime+}(x)=1$ for $x\geq b^{\ast}$ )(by \eqref{V.0bE.x}). 
Combining \eqref{gb(y)-gb(b).1} and \eqref{gb(y)-gb(b).2} we obtain
\begin{eqnarray}
\lim_{k\rightarrow\infty}\frac{g_{b^{\ast}}(y_{x,n_k}^+)-g_{b^{\ast}}(x)}{y_{x,n_k}^+-x}=0.\label{22823-6}
\end{eqnarray}
Since $y_{x,n_k}^+$ is a sub-sequence of $y_{x,n}^+$, we can infer from \eqref{22823-6} that \eqref{22823-4} holds. 

 From \eqref{v0b(x)-1} and \eqref{9823-8} we know $V_{0,b^{\ast}}^\prime(0)=\phi$ and $V_{0,b^{\ast}}^\prime(b^\ast)=1$. Then by  concavity of $V_{0,b^{\ast}}(x)$ on $(0,b^{\ast})$, we obtain
\begin{eqnarray}\label{1<V.pri<phi}
1\leq V_{0,b^{\ast}}^\prime(x)\leq\phi,\quad x\in(-\infty,\infty).
\end{eqnarray}
Hence, by noting that $R(s)$ is non-decreasing we obtain
\begin{align}
&\mathrm{E}_x\left[\int_{0-}^{t} e^{-(\delta+\lambda) s}V_{0,b^{\ast}}^{\prime}( {U}^{R,L}{(s-)}){d}R^c{(s)}\right]\nonumber\\
&+\mathrm{E}_x\left[\sum_{s\leq t}e^{-(\delta+\lambda) s}\big[V_{0,b^{\ast}}( {U}^{R,L}{(s-)}+\Delta  {}R{(s)})-V_{0,b^{\ast}}( {U}^{R,L}{(s-)})
\big]\mathbf{1}_{\{\Delta R{(s)}>0\}}\right]
\nonumber\\
\le  &\mathrm{E}_x\left[\int_{0-}^{t} e^{-(\delta+\lambda) s}\phi{d}R^c{(s)}\right]+\mathrm{E}_x\left[\sum_{s\leq t}e^{-(\delta+\lambda) s}\phi \Delta R_s\right]\nonumber\\
=&\phi \mathrm{E}_x\left[\int_{0-}^{t} e^{-(\delta+\lambda) s}{d}R^c{(s)}\right]. \label{23823-10}
\end{align}
Similarly, by noting $D(s)$ is non-decreasing we have
\begin{align}
&-\mathrm{E}_x\left[
\int_{0-}^{t} e^{-(\delta+\lambda) s}V_{0,b^{\ast}}^{\prime}( {U}^{R,L}{(s-)}){d}L^c{(s)})\right]\nonumber\\
&
+\mathrm{E}_x\left[\sum_{s\leq t}e^{-(\delta+\lambda) s}\big[V_{0,b^{\ast}}( {U}^{R,L}{(s-)}-\Delta  {}L{(s)})-V_{0,b^{\ast}}( {U}^{R,L}{(s-)})
\big]\mathbf{1}_{\{\Delta L{(s)}>0\}}\right]\nonumber\\
\le & -\mathrm{E}_x\left[
\int_{0-}^{t} e^{-(\delta+\lambda) s}{d}L^c{(s)}\right]
+\mathrm{E}_x\left[\sum_{s\leq t}e^{-(\delta+\lambda) s}(-\Delta L{(s)})\right]\nonumber\\
=&  -\mathrm{E}_x\left[
\int_{0-}^{t} e^{-(\delta+\lambda) s}{d}L{(s)}\right].
\label{23823-11}
\end{align}

By combining \eqref{23823-4}, \eqref{23823-7}, \eqref{g.b(x)<0}, \eqref{23823-10} and \eqref{23823-11} we can conclude that for any $t>0$ and $m\ge 1$,
\begin{align}
&\mathrm{E}_x\left[e^{-(\delta+\lambda)(t\wedge T_m)}V_{0,b^{\ast}}(U^{R,L}{(t\wedge T_m)})\right]-V_{0,b^{\ast}}(x)
\nonumber\\
\le &-\mathrm{E}_x\left[\int_{0-}^{t\wedge T_m} e^{-(\delta+\lambda) s}\lambda \beta V_{0,b^{\ast}}^E( {U}^{R,L}{(s-)}){d}s\right]+\mathrm{E}_x\left[\int_{0-}^{t\wedge T_m} e^{-(\delta+\lambda) s}\phi{d}R{(s)}\right]\nonumber\\
&-\mathrm{E}_x\left[\int_{0-}^{t\wedge T_m} e^{-(\delta+\lambda) s}{d}L{(s)})\right].
\label{23823-12}
\end{align}
That is, 
\begin{align}
V_{0,b^{\ast}}(x)\ge& \mathrm{E}_x\left[e^{-(\delta+\lambda)(t\wedge T_m)}V_{0,b^{\ast}}(U^{R,L}{(t\wedge T_m)})\right]+\mathrm{E}_x\left[\int_{0-}^{t\wedge T_m} e^{-(\delta+\lambda) s}\lambda \beta V_{0,b^{\ast}}^E( {U}^{R,L}{(s-)}){d}s\right]\nonumber\\
&-\mathrm{E}_x\left[\int_{0-}^{t\wedge T_m} e^{-(\delta+\lambda) s}\phi{d}R{(s)}\right]+\mathrm{E}_x\left[\int_{0-}^{t\wedge T_m} e^{-(\delta+\lambda) s}{d}L{(s)})\right].
\label{23823-13}
\end{align}
Note that under any admissible strategy $(R,L)$, $U^{R,L}({s})$ must be non-negative for $s>0$, and thus $U^{R,L}(s)\ge \min\{0,x\}$ under $\mathrm{P}_x$. As a result, by the non-decreasingness of $V_{0,b^\ast}$ it follows  that $V_{0,b^\ast}(U^{R,L}(s))\ge V_{0,b^\ast}(0\wedge x) $ under $\mathrm{P}_x$. Therefore,
\begin{align}
&\liminf_{t\rightarrow\infty}\liminf_{m\rightarrow\infty}\mathrm{E}_x\left[e^{-(\delta+\lambda)(t\wedge T_m)}V_{0,b^{\ast}}(U^{R,L}{(t\wedge T_m)})\right]\nonumber\\
\ge&  \liminf_{t\rightarrow\infty}\liminf_{m\rightarrow\infty}\mathrm{E}_x\left[e^{-(\delta+\lambda)(t\wedge T_m)}V_{0,b^{\ast}}(0\wedge x )\right]\nonumber\\
\ge&  \mathrm{E}_x\left[\liminf_{t\rightarrow\infty}\liminf_{m\rightarrow\infty}\left(e^{-(\delta+\lambda)(t\wedge T_m)}V_{0,b^{\ast}}(0\wedge x )\right)\right]=0,\label{24823-10}
    \end{align}
    where the last inequality follows by using the Fatou's Lemma twice. Note that 
\begin{align*}
&\mathrm{E}_x\left[\int_{0-}^{t\wedge T_m} e^{-(\delta+\lambda) s}\lambda \beta V_{0,b^{\ast}}^E( {U}^{R,L}{(s-)}){d}s\right]\nonumber\\
=&\mathrm{E}_x\left[\int_{0-}^{t\wedge T_m} e^{-(\delta+\lambda) s}\lambda \beta V_{0,b^{\ast}}^E( {U}^{R,L}{(s-)})\mathbf{1}_{\{V_{0,b^{\ast}}^E( {U}^{R,L}{(s-)})\ge 0\}}{d}s\right]\nonumber\\
&+\mathrm{E}_x\left[\int_{0-}^{t\wedge T_m} e^{-(\delta+\lambda) s}\lambda \beta V_{0,b^{\ast}}^E( {U}^{R,L}{(s-)}) \mathbf{1}_{\{V_{0,b^{\ast}}^E( {U}^{R,L}{(s-)})< 0\}}{d}s\right].
\end{align*}
By taking limits and then using the monotone convergence four times, we have
\begin{align}
&\liminf_{t\rightarrow\infty}\liminf_{m\rightarrow\infty}\mathrm{E}_x\left[\int_{0-}^{t\wedge T_m} e^{-(\delta+\lambda) s}\lambda \beta V_{0,b^{\ast}}^E( {U}^{R,L}{(s-)}){d}s\right]\nonumber\\
=&\mathrm{E}_x\left[\int_{0-}^{\infty} e^{-(\delta+\lambda) s}\lambda \beta V_{0,b^{\ast}}^E( {U}^{R,L}{(s-)}) \mathbf{1}_{\{V_{0,b^{\ast}}^E( {U}^{R,L}{(s-)})\ge 0\}}{d}s\right]\nonumber\\
&+\mathrm{E}_x\left[\int_{0-}^{\infty} e^{-(\delta+\lambda) s}\lambda \beta V_{0,b^{\ast}}^E( {U}^{R,L}{(s-)}) \mathbf{1}_{\{V_{0,b^{\ast}}^E( {U}^{R,L}{(s-)})< 0\}}{d}s\right]\nonumber\\
=&\mathrm{E}_x\left[\int_{0-}^{\infty} e^{-(\delta+\lambda) s}\lambda \beta V_{0,b^{\ast}}^E( {U}^{R,L}{(s-)}){d}s\right].\label{24823-11}
\end{align}
Recall that both $R$ and $D$ are non-decreasing. Then by applying the monotone convergence multiple times we obtain
\begin{align}
&\liminf_{t\rightarrow\infty}\liminf_{m\rightarrow\infty}\left( -\mathrm{E}_x\left[\int_{0-}^{t\wedge T_m} e^{-(\delta+\lambda) s}\phi{d}R{(s)}\right]+\mathrm{E}_x\left[\int_{0-}^{t\wedge T_m} e^{-(\delta+\lambda) s}{d}L{(s)}\right]\right)\nonumber\\
 =  &-\mathrm{E}_x\left[\int_{0-}^{\infty} e^{-(\delta+\lambda) s}\phi{d}R{(s)}\right]+\mathrm{E}_x\left[\int_{0-}^{\infty} e^{-(\delta+\lambda) s}{d}L{(s)}\right]. \label{24823-12}
    \end{align}
By combining \eqref{23823-13} and \eqref{24823-10}-\eqref{24823-12} we arrive at  
\begin{align}
V_{0,b^{\ast}}(x)\ge& \mathrm{E}_x\left[\int_{0-}^{\infty} e^{-(\delta+\lambda) s}\lambda \beta V_{0,b^{\ast}}^E( {U}^{R,L}{(s-)}){d}s\right]\nonumber\\
&-\mathrm{E}_x\left[\int_{0-}^{\infty} e^{-(\delta+\lambda) s}\phi{d}R{(s)}\right]+\mathrm{E}_x\left[\int_{0-}^{\infty} e^{-(\delta+\lambda) s}{d}L{(s)}\right]\nonumber\\
=&\mathcal{P}(x,\pi,\pi^{0,b^{\ast}}).
\end{align}
By the arbitrariness of $\pi=(R,L)\in\Pi$, we know that  
$V_{0,b^{\ast}}(x)\ge \sup_{\pi\in\Pi}\mathcal{P}(x,\pi,\pi^{0,b^{\ast}})$. This completes the  proof.
\hfill $\square$\\

\end{document}